\documentclass[final,10pt,journal,letterpaper]{IEEEtran}

\IEEEoverridecommandlockouts

\usepackage{amsmath}
\usepackage{amssymb}
\usepackage{amsfonts}
\usepackage{array}
\usepackage{dsfont}
\usepackage[compress]{cite}
\usepackage{url}
\usepackage{color}
\usepackage{graphicx}      
\usepackage{ntheorem}
\usepackage{tikz}
\usepackage{booktabs}

\makeatletter
\renewcommand*\env@matrix[1][*\c@MaxMatrixCols c]{%
  \hskip -\arraycolsep
  \let\@ifnextchar\new@ifnextchar
  \array{#1}}
\makeatother

\makeatletter
\g@addto@macro\normalsize{%
  \setlength\abovedisplayskip{4pt}
  \setlength\belowdisplayskip{4pt}
}
\makeatother

\theoremprework{%
\setlength\theorempreskipamount{0.5ex}\setlength\theorempostskipamount{0.5ex}
}
\newtheorem{definition}{Definition}
\theoremprework{%
\setlength\theorempreskipamount{1mm}\setlength\theorempostskipamount{1mm}
}

\theoremprework{%
\setlength\theorempreskipamount{0.5ex}\setlength\theorempostskipamount{0.5ex}
}

\theoremprework{%
\setlength\theorempreskipamount{1mm}\setlength\theorempostskipamount{1mm}
}

\theoremprework{%
\setlength\theorempreskipamount{0.5ex}\setlength\theorempostskipamount{0.5ex}
}
\newtheorem{proposition}{Proposition}
\theoremprework{%
\setlength\theorempreskipamount{1mm}\setlength\theorempostskipamount{1mm}
}

\theoremprework{%
\setlength\theorempreskipamount{0.5ex}\setlength\theorempostskipamount{0.5ex}
}
\newtheorem{lemma}{Lemma}
\theoremprework{%
\setlength\theorempreskipamount{1mm}\setlength\theorempostskipamount{1mm}
}

\theoremprework{%
\setlength\theorempreskipamount{0.5ex}\setlength\theorempostskipamount{0.5ex}
}
\newtheorem{remark}{Remark}
\theoremprework{%
\setlength\theorempreskipamount{1mm}\setlength\theorempostskipamount{1mm}
}

\theoremprework{%
\setlength\theorempreskipamount{0.5ex}\setlength\theorempostskipamount{0.5ex}
}
\newtheorem{assumption}{Assumption}
\theoremprework{%
\setlength\theorempreskipamount{1mm}\setlength\theorempostskipamount{1mm}
}

\theoremprework{%
\setlength\theorempreskipamount{0.5ex}\setlength\theorempostskipamount{0.5ex}
}
\newtheorem{theorem}{Theorem}
\theoremprework{%
\setlength\theorempreskipamount{1mm}\setlength\theorempostskipamount{1mm}
}

\theoremprework{%
\setlength\theorempreskipamount{0.5ex}\setlength\theorempostskipamount{0.5ex}
}

\theoremprework{%
\setlength\theorempreskipamount{1mm}\setlength\theorempostskipamount{1mm}
}

\theoremprework{%
\setlength\theorempreskipamount{0.5ex}\setlength\theorempostskipamount{0.5ex}
}

\theoremprework{%
\setlength\theorempreskipamount{1mm}\setlength\theorempostskipamount{1mm}
}

\newcommand{\revise}[1]{{\color{black}#1}}

\newcommand{\argmin}{\mathop{\text{arg\,min}}}

\newcommand{\Real}{{\mathbb R}} 
\newcommand{\Nat}{{\mathbb N}} 

\DeclareMathOperator{\lcm}{lcm}
\DeclareMathOperator{\modd}{mod}
\DeclareMathOperator{\dist}{dist}

\DeclareMathOperator{\proj}{proj}

\title{Secure and Private Implementation of Dynamic Controllers Using Semi-Homomorphic Encryption}

\author{Carlos Murguia\thanks{C.~Murguia and I.~Shames are with the Department of Electrical and Electronic Engineering at the University of Melbourne.}, Farhad Farokhi, \thanks{F.~Farokhi is with the CSIRO's Data61 and the Department of Electrical and Electronic Engineering at the University of Melbourne.}and Iman Shames}

\begin{document}
\maketitle

\begin{abstract} This paper presents a secure and private implementation of linear time-invariant dynamic controllers using Paillier's encryption, a semi-homomorphic encryption method. To avoid overflow or underflow within the encryption domain, the state of the controller is reset periodically. A control design approach is presented to ensure stability and optimize performance of the closed-loop system with encrypted controller.
\end{abstract}

\section{Introduction}
Internet of Things (IoT) has brought opportunities for flexibility of deployment and efficiency improvements. However, it threatens security and privacy of individuals and businesses as IoT devices, by design, share their information for processing over the cloud. This information can be secured from adversaries over the network by using encrypted communication channels~\cite{Patel2009ICS15387881538820}. This approach, although effective and necessary, does not address vulnerability of the data on servers running cloud-computing services. These services themselves can use the data for targeted advertisement or can be hacked for malicious purposes. Therefore, there is a need for a more secure methodology that addresses the security and privacy of data while being processed.

Thankfully secure cloud computing is possible with the use of homomorphic encryption methods -- encryption methods that allow computation over plain data by performing appropriate computations on the encrypted data~\cite{Gentry2009,Paillier1999,1057074}. The use of homomorphic encryption allows a controller to be remotely realised without needing to openly sharing private and sensitive data (and consenting to its use in an unencrypted manner). This paper specifically discusses secure and private implementation of linear time-invariant dynamic controllers with the aid of the Paillier's encryption~\cite{Paillier1999}, a semi-homomorphic encryption method.

The use of homomorphic encryption for secure control has been studied previously~\cite{kogiso2015cyber,farokhi2017secure,farokhi2016secure,  darup2018towards,kim2016encrypting,darup2018encrypted,lin2018secure, alexandru2018cloud}. However, all these studies consider static controllers. This is because, when dealing with dynamical control laws (with an encrypted memory/state that must be maintained remotely), the number of bits required for representing the state of the controller grows linearly with the number of iterations. This renders the memory useless after a few iterations due to an overflow or an underflow (i.e., number of fractional bits required for representing a number becomes larger than the  number of fractional bits in the fixed-point number basis). In fact, using rough calculations, it can be seen that for a system with sampling time of 10 milliseconds, 16 bits quantized controller parameters and measurements, and within an encryption space of \revise{2048} bits\footnote{Encryption keys with the length of \revise{2048} bits is recommended by National Institute of Standards and Technology (NIST) for data \revise{over 2016-2030}; see \revise{\url{https://www.keylength.com/en/4/}.} }, the state of the controller becomes incorrect after roughly \revise{1.2} seconds due to an overflow or underflow. \revise{The unstabilizing effect of restricting the memory of controllers to finite rings is illustrated in Section~\ref{sec:example} for the key length of 2048 bits using a controller that can easily stabilize a batch chemical reactor in the absence of encryption.}

\revise{There are multiple ways to deal with this issue:
\begin{enumerate}
\item We should decrypt the state of the encrypted controller, project it into the desired set of fixed-point rational numbers, and encrypt it again. To avoid this issue, the encrypted state can be sent to a trust third-party (e.g., an IoT device) to be decrypted, rounded, encrypted, and transmitted back. This adds unnecessary communication overhead and overburdens the computational units of the IoT device. Furthermore, by decrypting the state, the risk of a privacy or security breach increases.
\item We should restrict the controller parameters so that the state of the dynamic controller remains within the set of fixed-point rational numbers. This approach was pursued in~\cite{8619600}. This makes the problem of designing the controller into a mixed-integer optimization problem, which can be computationally exhaustive. However, a robust control approach can be taken to ensure that converting non-integer controllers to integer ones does not ruin stability~\cite{8619600}.
\item We should reset the controller, i.e., the state of the controller is set to a publicly known number (e.g., zero) periodically. In this case, the controller must be redesigned to ensure stability/performance, which this paper shows to remain a tractable optimisation problem. This is the approach chosen by the current paper.
\end{enumerate}}
\revise{In this paper, we only focus on encrypting the outputs of the system and the state of the controller. This is because the parameters of the controller are often  not sensitive in practice. For instance, in autonomous vehicles, the location and velocity are sensitive as they reveal private information about the user, e.g., home/work address and travel habits, while the controller parameters are implicitly related to the dynamics of the vehicle.}


Resetting controllers have been previously studied in~\cite{1100479,bakkeheim2008lyapunov,clegg1958nonlinear, krishnan1974synthesis,beker2004fundamental,prieur2018analysis, guo2012stability}. However, the synthesis approach in this paper is more general than those studies and further it is designed to accommodate challenges associated with the implementation of dynamical controllers over the cipher space. \revise{Particularly, majority of existing work on reset controllers focus on state dependent triggers. Due to the nature of our problem, where the controller cannot access to the unencrypted  state, those results are not applicable. Along the same lines, since we always have to reset the controller to the same state regardless of the state of the plant, the existing results for switched systems seem to be not applicable. }

The rest of the paper is organized as follows.  Preliminary materials on homomorphic encryption are presented in Section~\ref{sec:prelim}. The design and implementation of the controller is discussed in Section~\ref{sec:implementation}. Finally, numerical results are presented in Section~\ref{sec:example} and the paper is concluded in Section~\ref{sec:conclusions}. All proofs are presented in the appendices to improve the overall presentation of the paper.

\section{Preliminary Material} \label{sec:prelim}
In this paper, a tuple $(\mathbb{P},\mathbb{C},\mathbb{K},\mathfrak{E},\mathfrak{D})$ denotes a public key encryption scheme, where $\mathbb{P}$ is the set of plaintexts, $\mathbb{C}$ is the set of ciphertexts, $\mathbb{K}$ is the set of keys, $\mathfrak{E}$ is the encryption algorithm, and $\mathfrak{D}$ is the decryption algorithm. Each $\kappa=(\kappa_p,\kappa_s)\in\mathbb{K}$ is composed of a public key $\kappa_p$ (which is shared with and used by everyone for encrypting plaintexts) and a private key $\kappa_s$ (which is maintained only by the trusted parties for decryption). Algorithms $\mathfrak{E}$ and $\mathfrak{D}$ are publicly known while the keys, which set the parameters of these algorithms, are generated and used in each case. The use of the term ``algorithm'', instead of mapping or function, is due to the presence of random\footnote{These random elements are replaced with pseudo-random ones when implementing encryption and decryption algorithms.} elements in the encryption procedure possibly resulting in one plaintext being mapped to multiple ciphertexts. A necessary requirement for the encryption scheme is to be invertible, i.e., $\mathfrak{D}(\mathfrak{E}(x,\kappa_p),\kappa_p,\kappa_s)=x$ for all $x\in\mathbb{P}$ given $\kappa=(\kappa_p,\kappa_s)\in\mathbb{K}$.

\begin{definition}[Homomorphic \revise{Property}] \label{def:homomorphic} Assume there exist operators $\circ$ and $\diamond$ such that $(\mathbb{P},\circ)$ and $(\mathbb{C},\diamond)$ form groups.  A public key encryption $(\mathbb{P},\mathbb{C},\mathbb{K},\mathfrak{E},\mathfrak{D})$ is called called homomorphic if $\revise{\mathfrak{D}(}\mathfrak{E}(x_1,\kappa_p)\diamond\mathfrak{E}(x_1,\kappa_p)\revise{,\kappa_p,\kappa_s)}=x_1\circ x_2$ for all $x_1,x_2\in\mathbb{P}$ and $\kappa\in\mathbb{K}$.
\end{definition}

\revise{Throughout this paper, $|\mathbb{A}|$ denotes the cardinality of any set $\mathbb{A}$. Further, we define the notation $\mathbb{Z}_q:=\{0,\dots,q-1\}=\{n\modd q:\forall n\in\mathbb{Z}\}$} for all positive integers $q\in\mathbb{N}$.  \revise{In this paper, we assume that $\mathbb{P}=\mathbb{Z}_{n_p}$ and $\mathbb{C}=\mathbb{Z}_{n_c}$} with $n_p=|\mathbb{P}|$ and $n_c=|\mathbb{C}|$. A public key encryption $(\mathbb{P},\mathbb{C},\mathbb{K},\mathfrak{E},\mathfrak{D})$ is additively homomorphic if there exists an operator $\diamond$ such that Definition~\ref{def:homomorphic} is satisfied when the operator $\circ$ is defined as $x_1\circ x_2:=(x_1+x_2)\,\modd\,n_p$ for all $x_1,x_2\in\mathbb{P}$. For additively homomorphic schemes, in this paper, the notation $\oplus$ is used to denote the equivalent operator in the ciphertext domain ($\diamond$ in the definition above). Similarly, a public key encryption is multiplicatively homomorphic if there exists an operator $\diamond$ such that Definition~\ref{def:homomorphic} is satisfied with $\circ$ defined as $x_1\circ x_2:=(x_1x_2)\,\modd\,n_p$ for all $x_1,x_2\in\mathbb{P}$. If a public key encryption is both additively and multiplicatively homomorphic, it is fully homomorphic but, if only one of these conditions is satisfied, it is semi-homomorphic. Homomorphism shows there exist operations over ciphertexts that can generate encrypted versions of sumed or multiplied plaintexts without the need of decrypting their corresponding cuphertexts. An example of additively homomorphic encryption scheme is the Paillier's encryption method~\cite{Paillier1999}. ElGamal is an example of multiplicatively homomorphic encryption schemes~\cite{1057074}. Recently, several fully homomorphic encryption methods have been also developed, see, e.g.,~\cite{Gentry2009}.

Now, we define semantic security\revise{, borrowed from~\cite{katz2014introduction}}. A key $\kappa=(\kappa_p,\kappa_s)\in\mathbb{K}$ is randomly generated. A probabilistic polynomial time-bounded adversary proposes $x_1,x_2\in\mathbb{P}$. The agent chooses \revise{$x$} at random  from $\{x_1,x_2\}$ with equal probability, encrypts $x$ according to $y=\mathfrak{E}(x,\kappa_p)$, and sends $y$ to the adversary (along with the public key $\kappa_p$). The adversary produces $x'$, which is an estimate of $x$ based on all the avialable information (everything except $\kappa_s$, i.e., $x_1$, $x_2$, $y$, $\mathfrak{E}$, $\mathfrak{D}$, $\kappa_p$). The adversary's advantage (in comparison to that of a pure random number generator) is given by $\mathrm{Adv}(|\mathbb{K}|):=|\mathbb{P}\{x=x'\}-1/2|$. The public key encryption $(\mathbb{P},\mathbb{C},\mathbb{K},\mathfrak{E},\mathfrak{D})$ is semantically secure (alternatively known as indistinguishability under chosen plaintext attack) if $\mathrm{Adv}$ is negligible\footnote{A function $f:\mathbb{N}\rightarrow\mathbb{R}_{\geq 0}$ is called negligible if, for any $c\in\mathbb{N}$, there exists $n_c\in\mathbb{N}$ such that $f(n)\leq 1/n^c$ for all $n\geq n_c$.}.

In this paper, the results are presented for the Paillier's encryption method. It is \revise{noteworthy} that the Paillier's encryption method is semantically secure under the \emph{Decisional Composite Residuosity Assumption}, i.e., it is ``hard'' to decide whether there exists $y\in\mathbb{Z}_{N^2}$ such that $x= y^N \modd N$ for  $N\in\mathbb{Z}$ and $x\in\mathbb{Z}_{N^2}$. More information regarding the assumption can be found in \cite{Paillier1999,yi2014homomorphic}. This can be used to establish the security of the proposed framework.

The Paillier's encryption scheme is as follows. First the public and private keys are generated. To do so, large prime numbers $p$ and $q$ are selected randomly and independently of each other such that $\gcd(pq,(1-p)(1-q))=1$, where $\gcd(a,b)$ refers to the greatest common divisor of integers $a$ and $b$. The public key (which is shared with all the parties and is used for encryption) is $\kappa_p=pq$. The private key (which is only available to the entity that needs to decrypt the data) is $\kappa_s=(\lambda,\mu)$ with $\lambda=\lcm(p-1,q-1)$ and $\mu=\lambda^{-1}\modd \kappa_p$, where $\lcm(a,b)$ is the least common multiple of integers $a$ and $b$. The ciphertext of plain message $x\in\mathbb{P}=\mathbb{Z}_{\kappa_p}$ is
$\mathfrak{E}(x,\kappa_p)=(\kappa_p+1)^xr^{\kappa_p}\modd \kappa_p^2,$
where $r$ is randomly selected with uniform probability from $\mathbb{Z}^*_{\kappa_p}:=\{x\in\mathbb{Z}_{\kappa_p}|\gcd(x,\kappa_p)=1\}$. Finally, to decrypt any ciphertext $c\in\mathbb{C}=\mathbb{Z}_{\kappa_p^2}$, $\mathfrak{D}(c,\kappa_p,\kappa_s)=(L(c^\lambda \modd \kappa_p^2)\mu)\modd \kappa_p,$ where $L(z)=(z-1)/\kappa_p$.

\begin{proposition}\emph{\cite{Paillier1999}} 1) For $r,r'\in\mathbb{Z}^*_{\kappa_p}$ and $t,t' \in\mathbb{P}$ such that $t+t'\in\mathbb{P}$, $\mathfrak{E}(t,\kappa_p)\mathfrak{E}(t',\kappa_p)\modd \kappa_p^2=\mathfrak{E}(t+t',\kappa_p)$; 2) For $r\in\mathbb{Z}^*_{\kappa_p}$ and $t,t' \in\mathbb{P}$ such that $tt'\in\mathbb{P}$, $\mathfrak{E}(t,\kappa_p)^{t'}\modd \kappa_p^2=\mathfrak{E}(t't,\kappa_p)$.
 \label{prop:cipher}
\end{proposition}


Proposition~\ref{prop:cipher} shows that the Paillier's encryption is a semi-homomorphic encryption scheme, i.e., algebraic manipulation of the plain data is possible without decryption using appropriate operations over the encrypted data. The Paillier's encryption is additively homomorphic with operator $\oplus$ being defined as $x_1\oplus x_2=(x_1x_2)\modd \kappa_p^2$ for all $x_1,x_2\in\mathbb{C}$. Note that the Paillier's method is not multiplicatively homomorphic as $t'$ in the identity $\mathfrak{E}(t,\kappa_p)^{t'}\modd \kappa_p^2=\mathfrak{E}(t't,\kappa_p)$ in Proposition~\ref{prop:cipher} is not encrypted. Define $\triangle$ such that $x_1\triangle x_2=x_1^{x_2}\modd N^2$ for all $x_1\in\mathbb{C}$ and $x_2\in\mathbb{P}$. Note that $\triangle$ is not an operator (in the mathematical sense) as its operands belong to two difference sets; it is just a mapping.

\section{Dynamic Controller Implementation}
\label{sec:implementation}
Consider the discrete-time linear time invariant system
\begin{align}\label{eqn:system}
\mathcal{P}:
\left\{
\begin{array}{rlc}
x[k+1]\hspace*{-.1in}&=Ax[k]+Bu[k],& x[0]=x_0,\\
y[k]\hspace*{-.1in}&=Cx[k],
\end{array}
\right.
\end{align}
with $k \in \mathbb{N}$, state $x[k]\in\mathbb{R}^{n_x}$, control input $u[k]\in\mathbb{R}^{n_u}$, and output $y[k]\in\mathbb{R}^{n_y}$. \revise{Many linear time-invariant systems cannot be stabilized by static output feedback controllers~\cite{SYRMOS1997125,CAO19981641,Bara}. Therefore, dynamic output feedback controllers have been used for decades to stabilize system using only output measurements, e.g., standard Kalman-filter (or Luenberger observer) based linear regulators \cite{Astrom}, and general dynamic output feedback controllers for quadratic performance \cite{Scherer_IEEE}.} System \eqref {eqn:system} is controlled by a dynamic output feedback controller of the form
\begin{align}\label{eqn:controller}
\hspace{-.04in}\mathcal{C}\hspace{-.04in}:\hspace{-.04in}
\left\{\hspace*{-.07in}
\begin{array}{rl}
\hspace{-.03in}x_c[k+1]\hspace*{-.14in}&=\hspace{-.04in}
\begin{cases}
\hspace{-.04in}A_cx_c[k]+B_cy[k], & \hspace{-.1in}(k\hspace{-.02in}+\hspace{-.02in}1)\hspace{-.02in}\modd \hspace{-.02in}T\hspace{-.04in}\neq\hspace{-.02in} 0, \\
\hspace{-.03in}0, & \hspace{-.1in}(k\hspace{-.02in}+\hspace{-.02in}1)\hspace{-.02in}\modd \hspace{-.02in}T\hspace{-.04in}=\hspace{-.02in} 0,
\end{cases}
\\[1.5em]
u[k]\hspace*{-.14in}&=\hspace{-.04in}C_cx_c[k]+D_cy[k],
\end{array}
\right.
\end{align}
with controller state $x_c[k]\in\mathbb{R}^{n_c}$. It is assumed that the state of the controller resets every $T$ time steps, i.e., $x_c[\ell T]=0$ for all $\ell\in\mathbb{N}$. \emph{This is because implementing encrypted controllers over an infinite horizon is impossible due to memory issues (by multiplication of fractional numbers, the number of bits required for representing fractional and integer parts grow)}. 
Combining the dynamics in~\eqref{eqn:system} and~\eqref{eqn:controller} results in the augmented system:
\begin{align}
z[k+1]\hspace{-.04in}=\hspace{-.04in}
\begin{cases}
F(\mathcal{P},\mathcal{C})z[k], & (k\hspace{-.02in}+\hspace{-.02in}1)\hspace{-.02in}\modd \hspace{-.02in}T\hspace{-.02in}\neq\hspace{-.02in} 0, \\
\begin{bmatrix}
I & 0\\
0 & 0
\end{bmatrix}
F(\mathcal{P},\mathcal{C})z[k], & (k\hspace{-.02in}+\hspace{-.02in}1)\hspace{-.02in}\modd \hspace{-.02in}T\hspace{-.02in}=\hspace{-.02in} 0,
\end{cases}
\end{align}
where $z[k]:=\begin{bmatrix}
x[k]^\top &
x_c[k]^\top
\end{bmatrix}^\top$ and
\begin{align}
F(\mathcal{P},\mathcal{C}):=
\begin{bmatrix}
A+BD_cC & BC_c \\
B_cC & A_c
\end{bmatrix}. \label{Synthesis5b}
\end{align}
The following theorem provides a sufficient condition for the asymptotic stability of the origin of \eqref{eqn:system} in feedback with the resetting controller \eqref{eqn:controller}.
\begin{theorem} \label{tho:1} The closed-loop dynamics \eqref{eqn:system}-\eqref{eqn:controller} is globally asymptotically stable if there exist \revise{$P\in\mathbb{R}^{(n_c+n_x)\times (n_c+n_x)}$}, $\varepsilon\in(0,1)$, $\mu\in [-1,0)$, $\delta\in[1,\infty)$, and $\epsilon\in(0,\infty)$ satisfying:\begin{subequations} \label{eqn:tho:1:cond}
\begin{align}
&P\succ \epsilon I,\label{eqn:tho:1:cond1}\\
&F(\mathcal{P},\mathcal{C})^\top P F(\mathcal{P},\mathcal{C})\preceq (1+\mu) P,\label{eqn:tho:1:cond2}\\
& F(\mathcal{P},\mathcal{C})^\top \begin{bmatrix}
I & 0 \\
0 & 0
\end{bmatrix} P\begin{bmatrix}
I & 0 \\
0 & 0
\end{bmatrix} F(\mathcal{P},\mathcal{C})
\preceq \delta P,\label{eqn:tho:1:cond3}\\
&\delta(1+\mu)^{T-1}<\varepsilon.\label{eqn:tho:1:cond4}
\end{align}
\end{subequations}
\end{theorem}

\begin{IEEEproof} See Appendix~\ref{proof:tho:1}.
\end{IEEEproof}

The following result provides a sufficient condition for the stabilizability of the system using the resetting controller.

\begin{proposition} \label{prop:stars} If $n_c\geq n_x$, $(A,B)$ is stabilizable, and $(A,C)$ is detectable, there exist $\mu = \mu^*\in [-1,0)$ and $\epsilon = \epsilon^*\in(0,\infty)$ such that~\eqref{eqn:tho:1:cond1} and~\eqref{eqn:tho:1:cond2} are satisfied.
\end{proposition}

\begin{IEEEproof} See Appendix~\ref{proof:prop:stars}.
\end{IEEEproof}

For $\mu^*\in [-1,0)$ and $\epsilon^*\in(0,\infty)$ in Proposition~\ref{prop:stars}, the following problem can be solved to find the smallest resetting horizon $T$ for the dynamical controller:
\begin{subequations} \label{eqn:optimization}
\begin{align}
\min_{T\in\mathbb{N}}\hspace{-.1in} \min_{\scriptsize\begin{array}{c}
\varepsilon\in(0,1) \\ \delta\in[1,\infty)
\end{array}} & \hspace{-.1in}T,\\
\mathrm{s.t.}\hspace{.20in}&\hspace{-.1in}P\succ \epsilon^* I,\quad \delta(1+\mu^*)^{T-1}<\varepsilon,\\
&\hspace{-.1in}F(\mathcal{P},\mathcal{C})^\top PF(\mathcal{P},\mathcal{C})\hspace{-.04in}\preceq\hspace{-.04in} (1+\mu^*) P,\\
& \hspace{-.1in}F(\mathcal{P},\mathcal{C})^\top \hspace{-.04in}\begin{bmatrix}
I & 0 \\
0 & 0
\end{bmatrix}\hspace{-.04in}
P\hspace{-.04in}\begin{bmatrix}
I & 0 \\
0 & 0
\end{bmatrix}\hspace{-.04in} F(\mathcal{P},\mathcal{C})
\hspace{-.04in}\preceq \delta P.
\end{align}
\end{subequations}
\revise{Note that the conditions in Theorem~\ref{tho:1}, or optimization problem~\eqref{eqn:optimization},  are sufficient but not necessary. This is always the case when working with Lyapunov-based techniques for stability of dynamical systems \cite{Scherer_IEEE,Kha02}. In the next subsection, we provide change of variables to cast these conditions as linear matrix inequalities that can be solved off-line only once and passed to the cloud for real-time control.}

\subsection{Synthesis of Resetting Controllers}
\revise{In this subsection, we use appropriate change of variables to linearize the matrix inequalities in Theorem~\ref{tho:1} without generating conservatism.}
We provide tools for designing full order ($n_c=n_x$) resetting controllers of the form \eqref{eqn:controller} satisfying \eqref{eqn:tho:1:cond}. That is, we look for matrices $(A_c,B_c,C_c,D_c)$ satisfying the inequalities in \eqref{eqn:tho:1:cond} for some positive definite $P \in \Real^{2n_x \times 2n_x}$, $\mu \in [-1,0)$, $\delta \in (0,\infty)$, $\varepsilon \in (0,1)$, and $T \in \mathbb{N}$. Let $n_c = n_x$ and $P$ be positive definite. Consider $F(\mathcal{P},\mathcal{C})$ in \eqref{Synthesis5b} and define:
\begin{align}
&\begingroup
\renewcommand*{\arraycolsep}{3pt}
\tilde{F}(\mathcal{P},\mathcal{C}) :=  \begin{bmatrix} I & 0\\ 0 & 0 \end{bmatrix} F(\mathcal{P},\mathcal{C}) = \begin{bmatrix}  A+BD_cC  & BC_c \\ 0  & 0 \end{bmatrix}.\label{Synthesis5bb}\endgroup
\end{align}
For simplicity of notation, in this subsection, $F(\mathcal{P},\mathcal{C})$ and $\tilde{F}(\mathcal{P},\mathcal{C})$ are denoted by $F$ and $\tilde{F}$, respectively. Then, \eqref{eqn:tho:1:cond2} and \eqref{eqn:tho:1:cond3} can be written as
\begin{align}\label{Synthesis5c}
&F^\top  P F - (1+\mu)P \preceq 0,\quad \tilde{F}^\top  P \tilde{F} - \delta P \preceq 0,
\end{align}
where $0$ denotes the zero matrix of appropriate dimensions. Using properties of the Schur complement, inequalities \eqref{Synthesis5c} are fulfilled if and only if the following is satisfied:
\begin{align}\label{Synthesis5e}
&\mathcal{L} := \begin{bmatrix}  (1+\mu)P  & F^\top P \\ PF  & P \end{bmatrix} \succeq 0,\tilde{\mathcal{L}} := \begin{bmatrix}  \delta P  & \tilde{F}^\top P \\ P\tilde{F}  & P \end{bmatrix} \succeq 0.
\end{align}
Note that the blocks $PF$ and $P\tilde{F}$ are nonlinear functions of $(P,A_c,B_c,C_c,D_c)$. \revise{ In what follows, we propose a change of variables: $\left(P,A_c,B_c,C_c,D_c \right) \rightarrow \nu,$
so that, in the new variables $\nu$, we can obtain \emph{affine} matrix inequalities equivalent to \eqref{Synthesis5e}.} In particular, for positive definite $P$ and \emph{nonlinear} matrix inequalities $\mathcal{L} \geq0$ and $\tilde{\mathcal{L}}\geq 0$, we aim at finding two invertible matrices $\mathcal{Y}$ and $\mathcal{T}$, and variables $\nu$ such that the congruence transformations $P \rightarrow \mathcal{Y}^\top  P \mathcal{Y}$, $\mathcal{L} \rightarrow \mathcal{T}^\top \mathcal{L} \mathcal{T}$, and $\tilde{\mathcal{L}} \rightarrow \mathcal{T}^\top \tilde{\mathcal{L}} \mathcal{T}$ lead to new Linear Matrix Inequalities (LMIs) $\mathcal{Y}^\top  P \mathcal{Y} >0$, $\mathcal{T}^\top \mathcal{L} \mathcal{T} \geq 0$, and $\mathcal{T}^\top \tilde{\mathcal{L}} \mathcal{T} \geq 0$ in the variables $\nu$. Let $P$ be positive definite and partitioned as follows:
\begin{equation}\label{Synthesis1}
\begingroup
\renewcommand*{\arraycolsep}{1.5pt}
P := \begin{bmatrix}  X \hspace{2mm} & U  \\ U^\top   & \tilde{X} \end{bmatrix},\endgroup
\end{equation}
with $X,U,\tilde{X} \in \Real^{n_x \times n_x}$ and positive definite $X,\tilde{X}$. Define
\begin{equation}\label{Synthesis2}
\begingroup
\renewcommand*{\arraycolsep}{1.5pt}
\revise{P}^{-1} =: \begin{bmatrix}  Y \hspace{2mm} & V \\ V^\top   & \tilde{Y} \end{bmatrix}, \hspace{1mm} \mathcal{Y} := \begin{bmatrix}  Y \hspace{2mm} & I \\ V^\top   & 0 \end{bmatrix}, \hspace{1mm} \mathcal{Z} := \begin{bmatrix}  I  & 0 \\ X  & U \end{bmatrix}.\endgroup
\end{equation}
Using block matrix inversion formulas, it can be verified that $YX+VU^\top =I$ and $YU+V\tilde{X}=0$, which yields $\mathcal{Y}^\top  P = \mathcal{Z}$. Then,  $P \rightarrow \mathcal{Y}^\top  P \mathcal{Y}$ takes the form:
\begin{align}\label{Synthesis5}
&\begingroup
\renewcommand*{\arraycolsep}{3pt}
\mathcal{Y}^\top  P \mathcal{Y} = \mathcal{Z} \mathcal{Y} = \begin{bmatrix}  Y  & I \\ I  & X \end{bmatrix} =: \mathbf{P}(\nu).\endgroup
\end{align}
Define $\mathcal{T}:= \text{diag}[\mathcal{Y},\mathcal{Y}]$ with $\mathcal{Y}$ as introduced in \eqref{Synthesis2}. Then, the transformations $\mathcal{L} \rightarrow \mathcal{T}^\top \mathcal{L} \mathcal{T}$ and $\tilde{\mathcal{L}} \rightarrow \mathcal{T}^\top  \tilde{\mathcal{L}} \mathcal{T}$ can be written as
\begin{align}\label{Synthesis5fpp}
&\mathcal{T}^\top \mathcal{L} \mathcal{T} = \begin{bmatrix}  (1+\mu)\mathbf{P}(\nu)  & \mathcal{Y}^\top F^\top  \mathcal{Z}^\top  \\ \mathcal{Z}F\mathcal{Y}  & \mathbf{P}(\nu) \end{bmatrix},\\
&\mathcal{T}^\top  \tilde{\mathcal{L}} \mathcal{T} = \begin{bmatrix}  \delta \mathbf{P}(\nu)  & \mathcal{Y}^\top \tilde{F}^\top \mathcal{Z}^\top  \\ \mathcal{Z}\tilde{F}\mathcal{Y}  & \mathbf{P}(\nu) \end{bmatrix}. \label{Synthesis5fppp}
\end{align}
Using the structure of $F$ and $\tilde{F}$ and the change of variables:
\begin{align}
\label{change_of_coordinates}
&\begingroup
\hspace{-.1in}
\renewcommand*{\arraycolsep}{2pt}
\begin{pmatrix} K_1-XAY & K_2 \\ K_3 & K_4 \end{pmatrix}
\hspace{-.05in}:=\hspace{-.05in} \begin{pmatrix} U & XB \\ 0 & I_{n_u}  \end{pmatrix} \hspace{-.05in}\begin{pmatrix} A_c & B_c\\ C_c & D_c \end{pmatrix} \hspace{-.05in}\begin{pmatrix} V^\top  & 0 \\ CY & I_{n_y} \end{pmatrix} \hspace{-.05in}\endgroup,
\end{align}
the blocks $\mathcal{Z}F\mathcal{Y}$ and $\mathcal{Z}\tilde{F}\mathcal{Y}$ can be written as
\begin{align}
&\begingroup
\renewcommand*{\arraycolsep}{3pt}
\mathcal{Z}F\mathcal{Y} = \begin{bmatrix} AY + BK_3  & A+BK_4C \\ K_1 & XA + K_2C  \end{bmatrix} =: \mathbf{F}(\nu),\endgroup \label{Synthesis7}\\
&\begingroup
\renewcommand*{\arraycolsep}{3pt}
\mathcal{Z}\tilde{F}\mathcal{Y} = \begin{bmatrix} AY + BK_3  & A+BK_4C \\ XBK_3 + XAY  & XA + XBK_4C  \end{bmatrix} =: \tilde{\mathbf{F}}(\nu).\endgroup \label{Synthesis7B}
\end{align}
Therefore, under $\mathcal{T}$ and the change of variables in \eqref{change_of_coordinates}, we can write $\mathcal{T}^\top \mathcal{L}\mathcal{T}$ and $\mathcal{T}^\top  \tilde{\mathcal{L}}\mathcal{T}$ as follows:
\begin{align}\label{Synthesis5fpppp}
&\begingroup
\renewcommand*{\arraycolsep}{0pt}
\mathcal{T}^\top \mathcal{L} \mathcal{T} =  \begin{bmatrix}  (1+\mu)\mathbf{P}(\nu)  & \hspace{2mm} \mathbf{F}(\nu)^\top  \\ \mathbf{F}(\nu)  & \mathbf{P}(\nu) \end{bmatrix} =: \mathbf{L}(\nu), \endgroup \\[2mm]
&\begingroup
\renewcommand*{\arraycolsep}{0pt}
\mathcal{T}^\top  \tilde{\mathcal{L}} \mathcal{T} =  \begin{bmatrix}  \delta \mathbf{P}(\nu)  & \hspace{2mm} \tilde{\mathbf{F}}(\nu)^\top  \\ \tilde{\mathbf{F}}(\nu)  & \mathbf{P}(\nu) \end{bmatrix} =: \mathbf{S}(\nu), \endgroup \label{Synthesis5fpppppp}
\end{align}
with $\mathbf{P}(\nu),\mathbf{F}(\nu)$, and $\tilde{\mathbf{F}}(\nu)$ as defined in \eqref{Synthesis5}, \eqref{Synthesis7}, and \eqref{Synthesis7B}, respectively. Therefore, the original matrix inequality, $\mathcal{L} \succeq 0$ defined in \eqref{Synthesis5e}, that depends non-linearly on the decision variables $(P,A_c,B_c,C_c,D_c)$ is transformed into a new inequality, $\mathbf{L}(\nu) \succeq 0$, that is an affine function of the variables $\nu$. Note, however, that $\mathbf{S}(\nu) \succeq 0$ (the block $\tilde{\mathbf{F}}(\nu)$) is still nonlinear in the new variables $\nu$. In the following lemma, we give a sufficient condition, in terms of an affine inequality $\tilde{\mathbf{L}}(\nu) \succeq 0$, for $\mathbf{S}(\nu)$ to be positive semidefinite.

\begin{lemma}\label{lemma:bound}
Consider $\mathbf{P}(\nu)$ and $\mathbf{S}(\nu)$ defined in \eqref{Synthesis5} and \eqref{Synthesis5fpppppp}, respectively. Define the matrices:
\begin{align}
\mathbf{R}(\nu) &:= \begin{pmatrix} AY + BK_3 & A + BK_4C  \end{pmatrix}, \label{Synthesis31}\\[1mm]
\tilde{\mathbf{L}}(\nu) &:= \begin{bmatrix}  \delta \mathbf{P}(\nu) & \hspace{2mm} \mathbf{R}(\nu)^\top  \\ \mathbf{R}(\nu) & 2I_n - X \end{bmatrix}. \label{Synthesis34}
\end{align}
Then, $\tilde{\mathbf{L}}(\nu) \succeq 0  \Rightarrow \mathbf{S}(\nu) \succeq 0$.
\end{lemma}
\begin{IEEEproof} See Appendix~\ref{proof:lemma:bound}.
\end{IEEEproof}
Lemma \ref{lemma:bound} provides a sufficient condition, $\tilde{\mathbf{L}}(\nu) \succeq 0$, for the nonlinear matrix $\mathbf{S}(\nu)$ to be positive semidefinite. This $\tilde{\mathbf{L}}(\nu)$ is an affine function of $\nu$. Note, however, that finding $\nu$ satisfying $(\mathbf{P}(\nu) \succ 0,\hspace{.5mm} \mathbf{L}(\nu) \succeq 0,\hspace{.5mm}  \tilde{\mathbf{L}}(\nu) \succeq 0)$ might not be sufficient to guarantee the existence of $(P,A_c,B_c,C_c,D_c)$ satisfying $(P \succ 0,\hspace{.5mm} \mathcal{L} \succeq 0,\hspace{.5mm} \tilde{\mathcal{L}} \succeq 0)$. For this to be true, matrices $\mathcal{Y}$ and $\mathcal{T}$ must be invertible so that the transformations $P \rightarrow \mathcal{Y}^\top  P \mathcal{Y} = \mathbf{P}(\nu)$, $\mathcal{L} \rightarrow \mathcal{T}^\top \mathcal{L} \mathcal{T} = \mathbf{L}(\nu)$, and $\tilde{\mathcal{L}} \rightarrow \mathcal{T}^\top  \tilde{\mathcal{L}} \mathcal{T} = \mathbf{S}(\nu)$ are congruence transformations; and $\nu$ must render the change of variables in \eqref{change_of_coordinates} invertible.
\begin{lemma}\label{lemma:congruence}
Consider matrices $\mathcal{Y}$ and $\mathbf{P}(\nu)$ defined in \eqref{Synthesis2} and \eqref{Synthesis5}, respectively. Let $(X,Y)$ be such that $\mathbf{P}(\nu)\succ 0$. Then, $\mathcal{Y}$ and $\mathcal{T} = \text{diag}(\mathcal{Y},\mathcal{Y})$ are nonsingular and the change of variables in \eqref{change_of_coordinates} is invertible.
\end{lemma}
\begin{IEEEproof} See Appendix~\ref{proof:lemma:congruence}.
\end{IEEEproof}
Therefore, by Lemma \ref{lemma:congruence}, if $\mathbf{P}(\nu) \succ 0$, the transformations $P \rightarrow   \mathbf{P}(\nu)$, $\mathcal{L} \rightarrow  \mathbf{L}(\nu)$, and $\tilde{\mathcal{L}} \rightarrow  \mathbf{S}(\nu)$ are congruence transformations. The latter and the fact that (by Lemma \ref{lemma:bound}) $\tilde{\mathbf{L}}(\nu) \revise{\,\succeq\,} 0 \Rightarrow \mathbf{S}(\nu) \revise{\,\succeq\,} 0$ imply that
\begin{equation} \label{Synthesis11}
\left\{ \begin{array}{c}
( \mathbf{P}(\nu)  \revise{\,\succ\,} 0, \text{ }\mathbf{L}(\nu)  \revise{\,\succeq\,} 0,  \text{ and } \tilde{\mathbf{L}}(\nu)  \revise{\,\succeq\,} 0 )\\
\Downarrow \\
 ( P  \revise{\,\succ\,} 0, \text{ }\mathcal{L}  \revise{\,\succeq\,} 0,  \text{ and } \tilde{\mathcal{L}}  \revise{\,\succeq\,} 0 ),
\end{array} \right.
\end{equation}
for $P = \mathcal{Y}^{-\top}\mathbf{P}(\nu)\mathcal{Y}^{-1}$ and the controller matrices in \eqref{CONTROLLER} obtained by inverting \eqref{change_of_coordinates}. In the following lemma, we summarize the discussion presented above.

\begin{lemma}\label{synthesis_lemma1}
For given system matrices $(A,B,C)$. If there exist matrices $\nu = (X,Y,K_1,K_2,K_3,K_4)$, $K_2 \in \Real^{n_x \times n_u}$, $K_3 \in \Real^{n_y \times n_x}$, $K_4 \in \Real^{n_u \times n_y}$, $X,Y,K_1 \in \Real^{n_x \times n_x}$ satisfying $\mathbf{P}(\nu) \revise{\,\succ\,}0$, $\mathbf{L}(\nu) \revise{\,\succeq\,} 0$, and $\tilde{\mathbf{L}}(\nu) \revise{\,\succeq\,} 0$ with $\mathbf{P}(\nu)$, $\mathbf{L}(\nu)$, and $\tilde{\mathbf{L}}(\nu)$ as defined in \eqref{Synthesis5}, \eqref{Synthesis5fpppp}, and \eqref{Synthesis34}, respectively; then, there exist $(P,A_c,B_c,C_c,D_c)$ satisfying $P \revise{\,\succ\,}0$, $\mathcal{L}  \revise{\,\succeq\,} 0$, and $\tilde{\mathcal{L}}  \revise{\,\succeq\,} 0$ with $P$, $\mathcal{L}$, and $\tilde{\mathcal{L}}$ as defined in \eqref{Synthesis5e} and \eqref{Synthesis1}, respectively. Moreover, for every $\nu$ such that $\mathbf{P}(\nu)  \revise{\,\succ\,} 0$, $\mathbf{L}(\nu)  \revise{\,\succeq\,} 0$, and $\tilde{\mathbf{L}}(\nu)  \revise{\,\succeq\,} 0$, the change of variables in \eqref{change_of_coordinates} and matrix $\mathcal{Y}$ in \eqref{Synthesis2} are invertible and the $(P,A_c,B_c,C_c,D_c)$ obtained by inverting \eqref{Synthesis5} and \eqref{change_of_coordinates} are unique and satisfy the analysis inequalities \eqref{Synthesis5c}.
\end{lemma}

\begin{IEEEproof} See Appendix~\ref{proof:synthesis_lemma1}.
\end{IEEEproof}

By Lemma \ref{synthesis_lemma1}, the matrices $(P,A_c,B_c,C_c,D_c)$ obtained by inverting \eqref{Synthesis5} and \eqref{change_of_coordinates} satisfy inequalities \eqref{Synthesis5c} (and thus also \eqref{eqn:tho:1:cond2} and \eqref{eqn:tho:1:cond3}). Moreover, because the reconstructed $P$ is positive definite, inequality \eqref{eqn:tho:1:cond1} is satisfied with $\epsilon = \lambda_{\min}(P)$, where $\lambda_{\min}(P) \in \Real_{>0}$ denotes the smallest eigenvalue of $P$. Next, we give the synthesis result corresponding to Theorem~\ref{tho:1}.

\vspace{1mm}

\begin{theorem}\label{synthesis_theorem}
For given system matrices $(A,B,C)$ and constants $\varepsilon\in(0,1)$, $\mu\in [-1,0)$, $\delta\in[1,\infty)$, and $T \in \Nat$ satisfying \eqref{eqn:tho:1:cond4}, if there exist matrices $\nu = (X,Y,K_1,K_2,K_3,K_4)$ satisfying $\mathbf{P}(\nu)  \revise{\,\succ\,} 0$, $\mathbf{L}(\nu)  \revise{\,\succeq\,} 0$, and $\tilde{\mathbf{L}}(\nu)  \revise{\,\succeq\,} 0$; then, $P = \mathcal{Y}^{-\top}\mathbf{P(\nu)}\mathcal{Y}$ and the controller matrices in \eqref{CONTROLLER} satisfy the analysis inequalities \eqref{eqn:tho:1:cond1}-\eqref{eqn:tho:1:cond3} and thus render the closed-loop dynamics \eqref{eqn:system}-\eqref{eqn:controller} asymptotically stable.
\end{theorem}
\begin{IEEEproof} See Appendix~\ref{proof:synthesis_theorem}.
\end{IEEEproof}
\textbf{Controller Reconstruction.} For given $\nu$ satisfying the synthesis inequalities ($\mathbf{P}(\nu)  \revise{\,\succ\,} 0, \mathbf{L}(\nu)  \revise{\,\succeq\,} 0$, $\tilde{\mathbf{L}}(\nu)  \revise{\,\succeq\,} 0$):
\begin{enumerate}
  \item For given $X$ and $Y$, compute via singular value decomposition a full rank factorization $VU^\top  = I-YX$ with square and nonsingular $V$ and $U$.
  \item For given $\nu$ and invertible $V$ and $U$, solve the system of equations $\mathcal{Y}^\top  P \mathcal{T} = \mathbf{P}(\nu)$ and \eqref{change_of_coordinates} to obtain the matrices $(P,A_c,B_c,C_c,D_c)$.
\end{enumerate}

\begin{remark}
Note that $\varepsilon,\mu,\delta$, and $T$, in Theorem \ref{synthesis_theorem} must be fixed before looking for feasible solutions $\nu$ satisfying the synthesis LMIs: $\mathbf{P}(\nu)  \revise{\,\succ\,} 0$, $\mathbf{L}(\nu)  \revise{\,\succeq\,} 0$, and $\tilde{\mathbf{L}}(\nu)  \revise{\,\succeq\,} 0$. However, for any $\mu \in [-1,0)$ and $\delta\in[1,\infty)$, there always exist $\varepsilon \in (0,1)$ and $T \in \Nat$ satisfying \eqref{eqn:tho:1:cond4}. Moreover, the synthesis LMIs depend on $\nu$, $\delta$, and $\mu$ but not on $\varepsilon$ and $T$. Therefore, to find feasible controllers, we only have to fix $(\mu,\delta)$ and look for $\nu$ satisfying the synthesis LMIs. The constants $(\mu,\delta)$ are, in fact, variables of the synthesis problem; however, to linearize some of the constraints, we fix their value and search over $\mu \in [-1,0)$ and $\delta\in[1,\infty)$ to find feasible $\nu$. The latter increases the computations needed to find controllers; however, we can perform a bisection search over $\delta\in[1,\infty)$ and, because $\mu \in [-1,0)$ (a bounded set), a grid search over $\mu$ to decrease the required computations.
\end{remark}

\revise{Finally, note that the characteristics (e.g., unstable poles) of the system in~\eqref{eqn:system} make the feasibility of the design LMIs $\mathbf{P}(\nu)  \revise{\,\succ\,} 0$, $\mathbf{L}(\nu)  \revise{\,\succeq\,} 0$, and $\tilde{\mathbf{L}}(\nu)  \revise{\,\succeq\,} 0$ ``easier or harder'' for fixed resetting horizon $T$ and constants $\varepsilon,\mu,\delta$. Exploring this dependence, in general, is an avenue for future research.}

\subsection{Dynamic Controller Implementation}
\revise{In this subsection, we present the necessary transformations required for implementing encrypted dynamic control laws. } Before stating the next result, we introduce some notation. Define $\|A\|_{\max}:=\max_{i,j}|a_{ij}|$, where $a_{ij}$ denotes the entry in $i$-th row and $j$-th column of matrix $A$, and $\mathbb{Q}(n,m):=\{b\,|\,b=-b_n2^{n-m-1}+\sum_{i=1}^{n-1}2^{i-m-1}b_i,b_i\in\{0,1\}\,\forall i\in\{1,\dots,n\}\}.$
For any $x\in\mathbb{R}^q$ and $\mathbb{A}\subseteq\mathbb{R}^q$, let $ \proj(x,\mathbb{A})\in\argmin_{x'\in\mathbb{A}}\|x'-x\|_\infty$ and $\dist(x,\mathbb{A}):=\min_{x'\in\mathbb{A}}\|x'-x\|_\infty$. The quantization of $x\in\mathbb{R}^q$ is $\proj(x,\mathbb{Q}(n,m)^q)$ and the quantization error is $\|\proj(x,\mathbb{Q}(n,m)^q)-x\|_\infty=\dist(x,\mathbb{Q}(n,m)^q)$. The quantization of $X\in\mathbb{R}^{p\times q}$ is defined as $\proj(x,\mathbb{Q}(n,m)^{p\times q})\in\argmin_{x'\in\mathbb{A}}\|x'-x\|_{\max}$ and the quantization error as $\|\proj(x,\mathbb{Q}(n,m)^{p\times q})-x\|_{\max}$. Please refer to \cite{farokhi2017secure} for details about the quantization scheme.

\begin{theorem} \label{tho:2} Let 
\begin{subequations} \label{eqn:quantized}
\begin{align}
\bar{A}_c&=\proj(A_c,\mathbb{Q}(n,m)^{n_c\times n_c}),\\
\bar{B}_c&=\proj(B_c,\mathbb{Q}(n,m)^{n_c\times n_y}),\\
\bar{C}_c&=\proj(C_c,\mathbb{Q}(n,m)^{n_u\times n_x}),\\
\bar{D}_c&=\proj(D_c,\mathbb{Q}(n,m)^{n_u\times n_y}).
\end{align}
\end{subequations}
Then, there exists $\bar{n}\geq \bar{m}>0$ such that $F(\mathcal{P},\mathcal{C})$ satisfies~\eqref{eqn:tho:1:cond} if and only if $F(\mathcal{P},\bar{\mathcal{C}})$, where $\bar{\mathcal{C}}$ denotes the controller in~\eqref{eqn:controller} with quantized parameters $\bar{A}_c$, $\bar{B}_c$, $\bar{C}_c$, and $\bar{D}_c$ in~\eqref{eqn:quantized}, satisfies~\eqref{eqn:tho:1:cond} with the same $P$ for all $n\geq \bar{n}$ and $m\geq \bar{m}$.
\end{theorem}

\begin{IEEEproof} The proof follows from continuity of the eigenvalues. Note that, by increasing $n$ and $m$, the quantization error decreases (actually, it tends to zero).
\end{IEEEproof}

In what follows, we discuss the implementation of quantized resetting controllers using homomorphic encryption schemes and quantized sensor measurements. The controller, in this case, is given by
\begin{align} \label{eqn:controller_quantized}
\bar{\bar{\mathcal{C}}}\hspace{-.04in}:\hspace{-.04in}
\left\{\hspace*{-.07in}
\begin{array}{rl}\hspace{-.03in}
x_c[k+1]\hspace*{-.14in}&=\hspace{-.04in}
\begin{cases}
\hspace{-.04in}\bar{A}_cx_c[k]\hspace{-.02in}+\hspace{-.02in}\bar{B}_c\bar{y}[k], & \hspace{-.1in}(k\hspace{-.02in}+\hspace{-.02in}1)\hspace{-.02in}\modd \hspace{-.02in}T\hspace{-.03in}\neq\hspace{-.02in} 0, \\
\hspace{-.03in}0, & \hspace{-.1in}(k\hspace{-.02in}+\hspace{-.02in}1)\hspace{-.02in}\modd\hspace{-.02in} T\hspace{-.03in}=\hspace{-.02in} 0,
\end{cases}
\\[1.5em]
u[k]\hspace*{-.14in}&=\hspace{-.04in}\bar{C}_cx_c[k]+\bar{D}_c\bar{y}[k],
\end{array}
\right.
\end{align}
where $\bar{A}_c$, $\bar{B}_c$, $\bar{C}_c$, and $\bar{D}_c$ are defined in~\eqref{eqn:quantized} and
\begin{align}\label{eqn:output_quantized}
\bar{y}[k]&\in\argmin_{y\in \mathbb{Q}(n,m)^{n_y}} \|y-y[k]\|_{\infty}.
\end{align}
The difference between $\bar{\bar{\mathcal{C}}}$ in~\eqref{eqn:controller_quantized} and $\bar{\mathcal{C}}$ in Theorem~\ref{tho:2} is the quantization of the output measurements $y[k]$. The following standing assumption is made in this paper to ensure the stability of the closed-loop system.

\begin{assumption} \label{assum:1} $n\geq \bar{n}$ and $m\geq \bar{m}$ where $\bar{n}$ and $\bar{m}$ are given in Theorem~\ref{tho:2}.
\end{assumption}

The following theorem proves the stability of the system $\mathcal{P}$ with the quantized resetting controller $\bar{\bar{\mathcal{C}}}$. Note that, for any $r\in\mathbb{R}_{\geq 0}$, $\mathbb{B}(r):=\{x\,|\,\|x\|_2^2\leq r\}$.

\begin{theorem} \label{tho:3} Under Assumption~\ref{assum:1}, if there exist $\varepsilon\in(0,1)$, $\mu\in [-1,0)$, $\delta\in[1,\infty)$, $T \in \Nat$, and $\epsilon\in(0,\infty)$ such that inequalities in~\eqref{eqn:tho:1:cond} are satisfied and
\begin{align}\label{bound_on_n}
n>\log_2\bigg(\dfrac{\lambda_{\max}(C^\top C)}{\epsilon}x_0^\top \begin{bmatrix}
I & 0 \\
0 & 0
\end{bmatrix}P\begin{bmatrix}
I & 0 \\
0 & 0
\end{bmatrix}x_0\bigg)+1;
\end{align}
then, the system dynamics \eqref{eqn:system} with the quantized resetting controller in~\eqref{eqn:controller_quantized} is stable and, for some constant\footnote{See Appendix~\ref{proof:tho:3} for a description of this constant.} $\varrho>0$, $\lim_{k\rightarrow\infty} \dist(x[k],\mathbb{B}(\varrho2^{-m}))=0$.
\end{theorem}

\begin{IEEEproof} See Appendix~\ref{proof:tho:3}
\end{IEEEproof}

Theorem~\ref{tho:3} implies that the state of the system converges to a vicinity of the origin (instead of the origin itself) due to quantization effects. The volume of the this area can be arbitrarily reduced by increasing $m$ and thus the performance of the system can be arbitrarily improved.

\begin{lemma} \label{lemma:2} For the resetting quantized controller in~\eqref{eqn:controller_quantized}, $x_c[k]\in\mathbb{Q}((n_c+1)(k\modd T-1)+n_y+n(k\modd T+1),m(k\modd T+1))^{n_c}$ and $u_c[k]\in\mathbb{Q}((n_c+1)(k\modd T)+n_y+n(k\modd T+2),m(k\modd T+2))^{n_u}$.
\end{lemma}

\begin{IEEEproof} See Appendix~\ref{proof:lemma:2}.
\end{IEEEproof}

Using the change of variables:
\begin{subequations} \label{eqn:transformtointeger}
\begin{align}
\tilde{A}_c&=(2^m\bar{A}_c)\modd 2^{\tilde{n}},\\
\tilde{B}_c\revise{[k]}&=(\revise{2^{m(k\modd T+1)}}\bar{B}_c)\modd 2^{\tilde{n}},\\
\tilde{C}_c&=(2^m\bar{C}_c)\modd 2^{\tilde{n}},\\
\tilde{D}_c\revise{[k]}&=(\revise{2^{m(k\modd T+1)}}\bar{D}_c)\modd 2^{\tilde{n}},\\
\tilde{x}_c[k]&=(2^{m(k\modd T+1)}\bar{x}_c[k])\modd 2^{\tilde{n}},\\
\tilde{y}[k]&=(2^{m(k\modd T+1)}\bar{y}[k])\modd 2^{\tilde{n}},\\
\tilde{u}[k]&=(2^{m(k\modd T+2)}\bar{u}[k])\modd 2^{\tilde{n}},
\end{align}
\end{subequations}
with $\tilde{n}>(n_c+1)T+n_u+n(T+2)$, the resetting quantized controller in~\eqref{eqn:controller_quantized} can be rewritten as
\begin{align}  \label{eqn:controller_positiveint}
\hspace{-.1in}\tilde{\mathcal{C}}\hspace{-.04in}:\hspace{-.04in}
\left\{\hspace*{-.07in}
\begin{array}{rl}\hspace{-.04in}
\tilde{x}_c[k+1]\hspace*{-.14in}&=\hspace{-.04in}
\begin{cases}
\hspace{-.04in}\tilde{A}_c\tilde{x}_c[k]\hspace{-.02in}+\hspace{-.02in}\tilde{B}_c\revise{[k]}\tilde{y}[k], & \hspace{-.1in}(k\hspace{-.02in}+\hspace{-.02in}1)\hspace{-.02in}\modd \hspace{-.02in}T\hspace{-.03in}\neq\hspace{-.02in} 0, \\
\hspace{-.03in}0, & \hspace{-.1in}(k\hspace{-.02in}+\hspace{-.02in}1)\hspace{-.02in}\modd\hspace{-.02in} T\hspace{-.03in}=\hspace{-.02in} 0,
\end{cases}
\\[1.5em]
\tilde{u}[k]\hspace*{-.14in}&=\hspace{-.04in}\tilde{C}_c\tilde{x}_c[k]+\tilde{D}_c\revise{[k]}\tilde{y}[k].
\end{array}
\right.
\end{align}
Note that, by Lemma~\ref{lemma:2}, $\tilde{A}_c,\tilde{B}_c,\tilde{C}_c,\tilde{D}_c,\tilde{x}_c,\tilde{y},\tilde{u}$ are \textit{positive integers}. This is useful because the Paillier's scheme can only work with finite ring of positive integers. Therefore, the update equation can now be implemented using Paillier's encryption scheme. The correctness of this implementation follows from the results of~\cite{farokhi2017secure} on fixed-point rational numbers.

\begin{figure}[t]
  \centering
  \includegraphics[scale=.13]{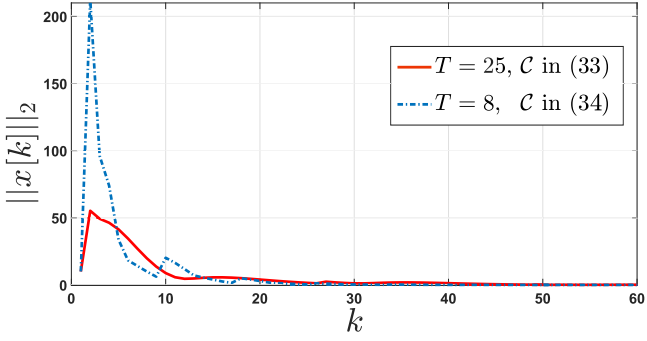}
  \vspace{-.1in}
  \caption{Norm of the state of the closed-loop system  $\|x[k]\|_2$ with quantized controller \eqref{eqn:controller_quantized} and quantizer resolution $(n,m) = (24,14)$.}\label{Fig1}
\end{figure}

\begin{figure}[t]
  \centering
  \includegraphics[scale=.13]{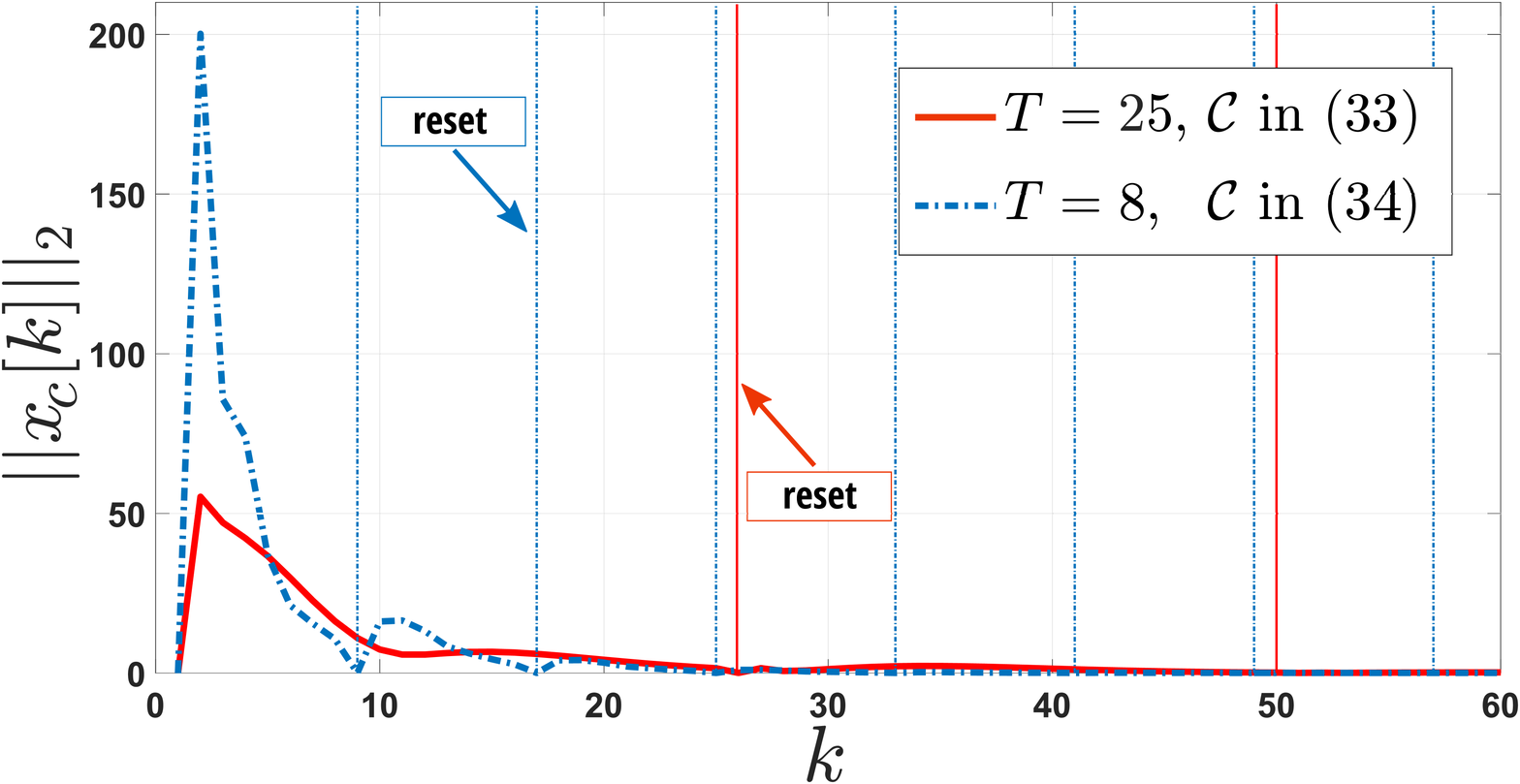}
  \vspace{-.1in}
  \caption{Norm of the state of the quantized controller $\|x_c[k]\|_2$  in \eqref{eqn:controller_quantized}  with quantizer resolution $(n,m) = (24,14)$.}\label{Fig2}
  \vspace{-.2in}
\end{figure}

First, the public and private keys must be generated such that $\kappa_p\geq 2^{\tilde{n}+1}$ to ensure that no unintended overflow occurs when using the encrypted numbers. The sensors measure, quantize, and encrypt the output to obtain
\begin{align}\label{encrypted_y}
\check{y}_i[k]:=\mathfrak{E}(\tilde{y}_i[k],\kappa_p).
\end{align}
The controller follows the encrypted version of~\eqref{eqn:controller_positiveint} to update its state  and compute the actuation signal as
\begin{align}\label{encrypted_xc}
(\check{x}_c)_i[k+1]=&
\begin{cases}
\bigg(\oplus_{j=1}^{n_x} (\check{x}_c)_j[k]\triangle(\tilde{A}_c)_{ij}\bigg)\\
\hspace{.2in}\oplus \bigg(\oplus_{j=1}^{n_y} \check{y}_j[k]\triangle(\tilde{B}_c)_{ij}\bigg),\\
&\hspace{-1in} (k+1)\modd T\neq 0,\\
\mathfrak{E}(0,\kappa_p), &\hspace{-1in} (k+1)\modd T=0,
\end{cases}
\end{align}
\begin{align}\label{encrypted_u}
\check{u}_i[k]=&\bigg(\oplus_{j=1}^{n_c} (\check{x}_c)_j[k]\triangle(\tilde{C}_c)_{ij}\bigg)\nonumber\\
&\oplus \bigg(\oplus_{j=1}^{n_y} (\check{y})_j[k]\triangle(\tilde{D}_c)_{ij}\bigg).
\end{align}
Finally, the actuator extract the control signal by
$\tilde{u}_i[k]=\mathfrak{D}(\check{u}_i[k],\kappa_p,\kappa_s)\modd 2^{\tilde{n}},$
and implements
$u_i[k]=2^{-m(k\modd T+2)}(\tilde{u}_i[k]-2^{\tilde{n}}\mathds{1}_{\tilde{u}_i[k]\geq 2^{\tilde{n}-1}}).$

\revise{
\begin{remark}
National Institute of Standards and Technology (NIST) recommends the use of key length of 2048 bits for factoring-based asymmetric encryption to guarantee that brute-force attacks are not physically possible during the life-time of the services based on projections of computing technologies.  This high standard might not be necessary for some applications, such as remote control of autonomous vehicles. To demonstrate this, consider RSA,  which is a similar encryption methodology and also a semi-homomorphic encryption relying on hardness of prime number factorization. RSA encryption has been attacked repeatedly using a brute-force methodology; see RSA Challenge\footnote{\url{https://en.wikipedia.org/wiki/RSA_Factoring_Challenge}}. Factorization of 663 bit numbers has been shown to take approximately 55 CPU-Years\footnote{A CPU-Year is the amuont of computing work done by a 1 Giga Floating Point Operations Per Second (FLOP) reference machine in a year of dedicated service (8760 hours).}~\cite{elbirt2009understanding}. Using IBM Watson (used recently for natural language processing to win quiz show Jeopardy), factorization of 663 bit numbers  takes approximately 2.5 years. These numbers are certainly not safe for use in finance or military. However, for remote control of autonomous vehicles, these keys may provide strong-enough guarantees as, by the time that an adversary breaks the code, the autonomous vehicle is in an entirely different location.
\end{remark}
 }

\begin{figure}[t]
  \centering
  \includegraphics[scale=.1275]{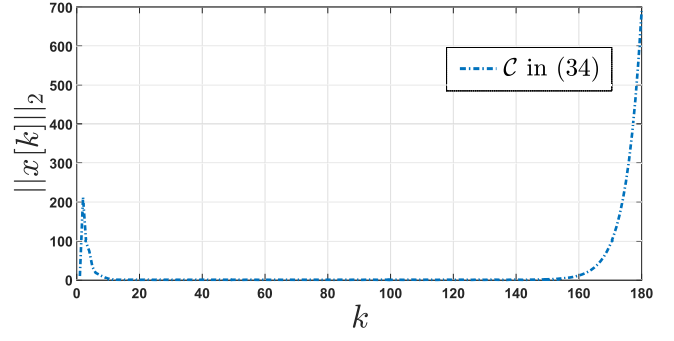}
  \vspace{-.1in}
  \caption{Norm of the state of the closed-loop system  $\|x[k]\|_2$ with quantized controller \eqref{eqn:controller_quantized2} and quantizer resolution $(n,m) = (24,14)$.}\label{Fig3}
\end{figure}

\begin{figure}[t]
  \centering
  \includegraphics[scale=.1275]{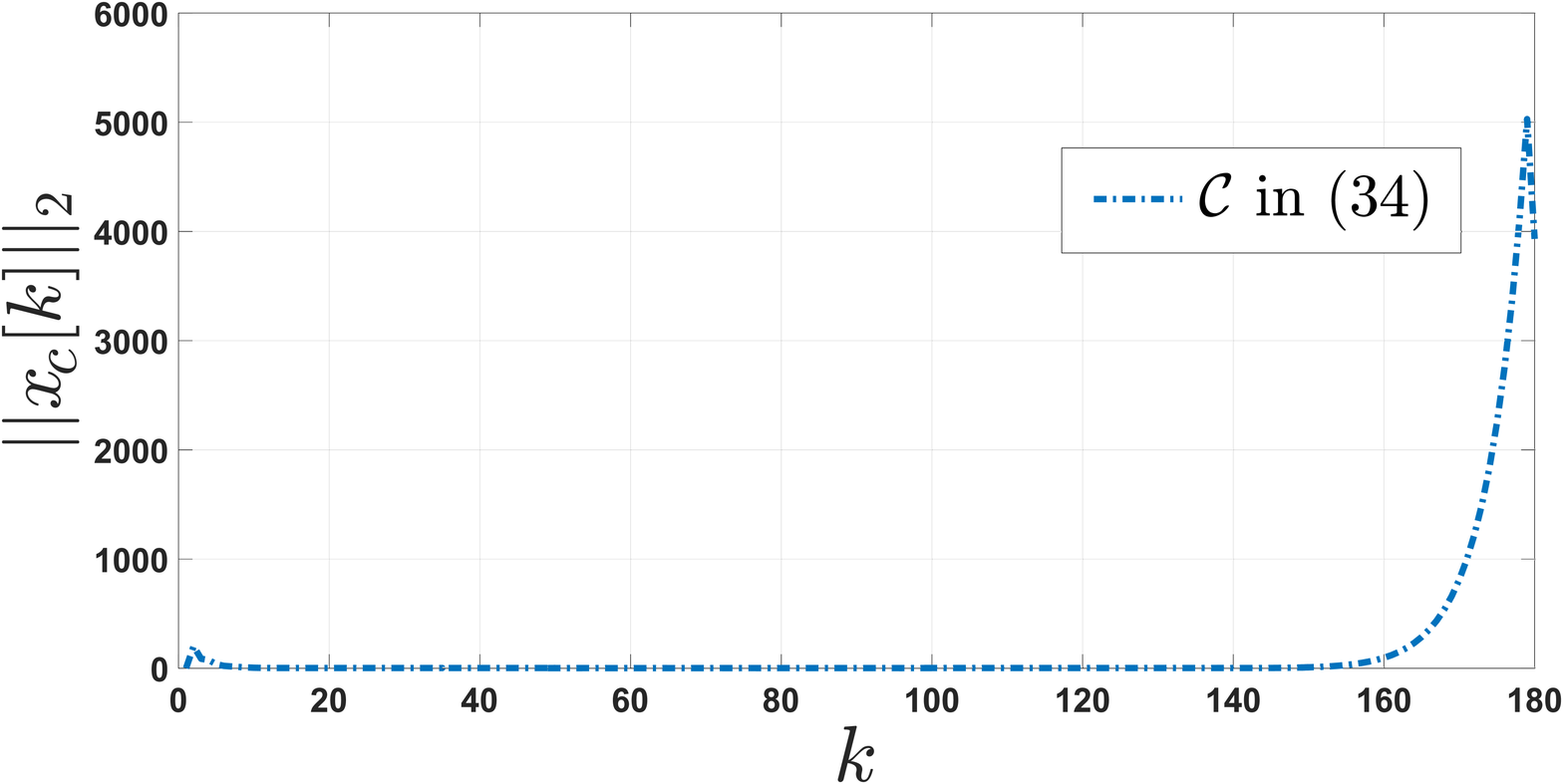}
  \vspace{-.1in}
  \caption{Norm of the state of the quantized controller $\|x_c[k]\|_2$  in \eqref{eqn:controller_quantized2}  with quantizer resolution $(n,m) = (24,14)$.}\label{Fig4}
  \vspace{-.2in}
\end{figure}

\section{Case Study of a Chemical Batch Reactor} \label{sec:example}

\revise{We illustrate the performance of our results through a case study of a batch chemical reactor. This case study has been developed over the years as a benchmark example for networked control systems, see, e.g., \cite{NCS1,NCS3,NCS2}. The reactor considered here is open-loop unstable, has one input, and two outputs (please refer to \cite{NCS1,NCS3,NCS2} for details about the system dynamics). We exactly discretize the reactor dynamics introduced in \cite{NCS3} with sampling period $h=0.1$. The resulting discrete-time linear system is of the form \eqref{eqn:system} with matrices $A,B,C$ as follows:
\begin{align}\label{Simul1b}
&\hspace{-.1in}\begin{bmatrix}[c|c|c]
  A & B & C^\top
\end{bmatrix}\nonumber\\
&= \hspace{-.05in}
\begin{bmatrix}[cccc|c|cc]
    \hspace{2.5mm}1.18  &  \hspace{2.5mm}0    &   \hspace{2.5mm}0.51 &  -0.40  &  0 & \hspace{2.5mm}1 & \hspace{.5mm}0\\
   -0.05  &  \hspace{2.5mm}0.66 &  -0.01 &   \hspace{2.5mm}0.06  &  0.47 & \hspace{2.5mm}0 & \hspace{.5mm}1\\
    \hspace{2.5mm}0.08  &  \hspace{2.5mm}0.34 &   \hspace{2.5mm}0.56 &   \hspace{2.5mm}0.38  &  0.21 & \hspace{2.5mm}1 & \hspace{.5mm}0\\
    \hspace{2.5mm}0     &  \hspace{2.5mm}0.34 &   \hspace{2.5mm}0.09 &   \hspace{2.5mm}0.85  &  0.21 & -1 & \hspace{.5mm}0
\end{bmatrix} \hspace{-.03in}.
\end{align}
Note that $\text{ eig}[A] = \{1.22,1.01,.60,.42  \}$; thus the system is open-loop unstable. Moreover, it can be verified (e.g., using the tools in \cite[Theorem 3.3]{Bara}) that there does not exist a static output feedback controllers of the form $u[k] = Ly[k]$, $L \in \mathbb{R}^{1 \times 2}$, stabilizing system \eqref{eqn:system} with matrices $(A,B,C)$ as in \eqref{Simul1b}. First, using the synthesis results in Section \ref{sec:implementation}, we design switching dynamic output feedback controllers of the form \eqref{eqn:controller}. Using Theorem \ref{synthesis_theorem}, and conducting a bisection in $\delta \in [1,\infty)$, and a line search in $\mu \in [-1,0)$, we look for the smallest $\delta$ for which there exist $\mu \in [-1,0)$ and $\nu$ satisfying the synthesis LMIs in Theorem \ref{synthesis_theorem}. The obtained $\delta$ is given by $\delta= \delta^* = 55.0$, the corresponding $\mu$ is $\mu = \mu^* = -0.15$, the resetting horizon is $T^* = \text{argmin}_{T \in \Nat} \hspace{1mm} \delta^*(1+\mu^*)^T = 25$, and the reconstructed $A_c,B_c,C_c,D_c$ (see Section \ref{sec:implementation}) are given in
\begin{align}
 &\begin{bmatrix}[c|c]
  A_c & B_c\\ \cmidrule(lr){1-2}
  C_c & D_c
\end{bmatrix}\nonumber \\
&\hspace{-.05in}=\hspace{-.05in}
\begin{bmatrix}[cccc|cc]
    \hspace{2.5mm}0.26 &  -0.03 &  -0.29 &   \hspace{2.5mm}0.31 &  -0.52 &  -0.03\\
   -0.32 &   \hspace{2.5mm}1.24 &   \hspace{2.5mm}1.40 &  -3.05 &   \hspace{2.5mm}5.46 &   \hspace{2.5mm}1.25\\
   -0.45 &   \hspace{2.5mm}0.02 &   \hspace{2.5mm}0.87 &  -0.75 &   \hspace{2.5mm}2.32 &  -0.01\\
   -0.05 &  -0.04 &   \hspace{2.5mm}0.72 &  -0.51 &   \hspace{2.5mm}2.28 &  -0.08\\ \cmidrule(lr){1-6}
    \hspace{2.5mm}1.02 &  -2.65 &  -2.65 &   \hspace{2.5mm}6.28 & -11.3 &  -4.09
\end{bmatrix}\hspace{-.05in}.	\label{Simul2b}
\end{align}
This controller satisfies the original inequalities in \eqref{eqn:tho:1:cond} with $\epsilon = 0.0026$ and any $\varepsilon \in (0.9459,1)$. For comparison, let $\mu = \mu^* = -0.65$ and search for the smallest $\delta$ for which there exists $\nu$ satisfying the synthesis LMIs in Theorem \ref{synthesis_theorem}. In this case, $\delta^* = 3000$, the smallest resetting horizon is $T^* = \text{argmin}_{T \in \Nat} \hspace{1mm} \delta^*(1+\mu^*)^T = 8$, the reconstructed $A_c,B_c,C_c,D_c$ are in
\begin{align}
 &\begin{bmatrix}[c|c]
  A_c & B_c\\ \cmidrule(lr){1-2}
  C_c & D_c
\end{bmatrix}\nonumber\\
& \hspace{-.05in}=\hspace{-.05in}
\begin{bmatrix}[cccc|cc]
   -0.18 &  -0.01 &  -0.77 &   \hspace{2.5mm}0.84 &  -1.11 &  -0.01\\
    \hspace{2.5mm}9.17 &   \hspace{2.5mm}0.43 &   \hspace{2.5mm}13.4 &  -16.2 &   \hspace{2.5mm}22.8 &   \hspace{2.5mm}0.42\\
    \hspace{2.5mm}1.24 &   \hspace{2.5mm}0.10 &   \hspace{2.5mm}3.82 &  -4.22 &   \hspace{2.5mm}7.81 &   \hspace{2.5mm}0.06\\
    \hspace{2.5mm}1.32 &   \hspace{2.5mm}0.10 &   \hspace{2.5mm}3.47 &  -3.87 &   \hspace{2.5mm}7.89 &   \hspace{2.5mm}0.06\\ \cmidrule(lr){1-6}
   -19.6 &  -0.93 &  -28.8 &   \hspace{2.5mm}34.9 &  -49.0 &  -2.33
\end{bmatrix}\hspace{-.05in}.\label{Simul3b}
\end{align}
This controller satisfies the original inequalities in \eqref{eqn:tho:1:cond} with
 $\epsilon =2.8 \times 10^{-5}$, and $\varepsilon \in (0.6756,1)$.

Next, we quantize the controller matrices according to \eqref{eqn:quantized} to obtain $(\bar{A}_c,\bar{B}_c,\bar{C}_c,\bar{D}_c)$ with quantizer resolution $(n,m)=(24,14)$. It can be verified that for $A_c,B_c,C_c,D_c$ in \eqref{Simul2b} and \eqref{Simul3b}, the corresponding $\bar{A}_c,\bar{B}_c,\bar{C}_c,\bar{D}_c$ satisfy the conditions of Theorem~\ref{tho:3} with $(n,m)=(24,14)$. We quantize sensor measurements $y[k]$ according to \eqref{eqn:output_quantized} with the same resolution $(n,m)=(24,14)$, and close the system dynamics with the quantized controller in \eqref{eqn:controller_quantized}. By Theorem~\ref{tho:3}, the quantizer resolution must satisfy inequality \eqref{bound_on_n} to ensure practical stability of \eqref{eqn:system} in feedback with \eqref{eqn:controller_quantized} in the sense of Theorem~\ref{tho:3}. Inequality \eqref{bound_on_n}, with initial condition $[x(0)^T,x_c(0)^T] = [-6.83,-5.18,-4.05,-3.12,0,0,0,0]$, amounts to $n>17$ for the controller in \eqref{Simul2b} and to $n>23$ for the controller in \eqref{Simul3b}. Therefore, $(n,m)=(24,14)$ is enough for practical stabilization using the controllers in \eqref{Simul2b} and \eqref{Simul3b}. Figures~\ref{Fig1} and~\ref{Fig2} show $\|x(k)\|_2$ and $\|x_c(k)\|_2$ of the closed-loop dynamics for quantized controllers corresponding to the controllers in \eqref{Simul2b} and \eqref{Simul3b} with $(n,m)=(24,14)$.

To illustrate the need for the proposed resetting controller, we naively implement a standard quantized dynamic controller of the form:
\begin{align} \label{eqn:controller_quantized2}
\left\{
\begin{array}{l}
\bar{x}_c[k+1] = \bar{A}_c\bar{x}_c[k] + \bar{B}_c\bar{y}[k],\\[1mm]
\hspace{7.5mm}\bar{u}[k] = \bar{C}_c\bar{x}_c[k]+\bar{D}_c\bar{y}[k].
\end{array}
\right.
\end{align}
We use the same quantizer resolution $(n,m)=(24,14)$, and compute the matrices $(\bar{A}_c,\bar{B}_c,\bar{C}_c,\bar{D}_c)$ using \eqref{eqn:quantized} with $(A_c,B_c,C_c,D_c)$ from \eqref{Simul3b}. This controller is stabilizing even without resets. Note that the Paillier's encryption only works over the ring of positive integers $\mathbb{Z}_{\kappa_p}$, and thus the controller needs to be transformed so that its states and parameters always belong to this ring. Therefore, as also required for the resetting controller, we must transform $\bar{y}[k]$ and $(\bar{A}_c,\bar{B}_c,\bar{C}_c,\bar{D}_c)$ into positive integers, which can be done using the change of variables in~\eqref{eqn:transformtointeger} replacing $k\modd T$ with $k$ (as there is no resetting after $T$ steps in this case). Let the integer representations of $\bar{y}[k]$ and $(\bar{A}_c,\bar{B}_c,\bar{C}_c,\bar{D}_c)$ be similarly denoted by $\tilde{y}[k]$ and $(\tilde{A}_c,\tilde{B}_c,\tilde{C}_c,\tilde{D}_c)$. The equivalent controller in the integer domain is then given by
\begin{align} \label{eqn:controller_integer2}
\left\{
\begin{array}{l}
\tilde{x}_c[k+1] = \big( \tilde{A}_c\tilde{x}_c[k] + \tilde{B}_c\tilde{y}[k] \big)\modd 2^{\tilde{n}},\\[1mm]
\hspace{7.5mm}\tilde{u}[k] = \big( \tilde{C}_c\tilde{x}_c[k]+\tilde{D}_c\tilde{y}[k] \big)\modd 2^{\tilde{n}}.
\end{array}
\right.
\end{align}
Finally, given $\tilde{u}[k]$, the actuators implement the control action:
\begin{align} \label{eqn:controller_integer2_extractred}
u_i[k]= 2^{-m(k+2)}(\tilde{u}_i[k]-2^{\tilde{n}}\mathds{1}_{\tilde{u}_i[k]\geq 2^{\tilde{n}-1}}).
\end{align}
Because we must ensure that $2^{\tilde{n}}\leq \kappa_p$, we need to select a large, yet finite $\tilde{n}$. Here, for illustration purposes, we selected a key length of 2048\,bits and $\tilde{n} = 2014$. Figure~\ref{Fig3} and~\ref{Fig4} illustrate the norm of the state of the closed-loop system, $\|x[k]\|_2$, with controller \eqref{Simul3b}-\eqref{eqn:controller_integer2_extractred}, and the state of the controller in the quantized domain, $\|x_c[k]\|$, respectively. Note that, even though \eqref{Simul3b}-\eqref{eqn:controller_quantized2} is a stabilizing controller (if no under/over flow occur), when implementing \eqref{Simul3b}-\eqref{eqn:controller_integer2_extractred}, the closed-loop system is unstable due to under/over flows. 
}

\section{Conclusions and Future Work} \label{sec:conclusions}
A secure and private implementation of linear time-invariant dynamic controllers using the Paillier's encryption was presented. The state is reset to zero periodically to avoid data overflow or underflow within the encryption space. A control design approach was presented to ensure the stability and performance of the closed-loop system with encrypted controller. Future work can focus on nonlinear dynamical systems and controllers.

\appendices

\section{Proof of Theorem~\ref{tho:1}}\label{proof:tho:1}
For the sake brevity $F(\mathcal{P},\mathcal{C})$ is denoted by $F$. Consider the quadratic Lyapunov function $V(z[k])=z[k]^\top Pz[k]$, where $P$ is a positive semi-definite matrix. For all $k\in\{\ell T,\dots,(\ell+1)T-2\}$ with $\ell\in\mathbb{N}$ the Lyapunov function evolves as
\begin{align}
V(z[k+1])-V(z[k])&=z[k]^\top (F^\top PF-P)z[k]\nonumber\\
&\leq \mu z[k]^\top P z[k]\nonumber\\
&= \mu V(z[k]). \label{eqn:proof:0}
\end{align}
Therefore,
\begin{align} \label{eqn:proof:1}
V(z[(\ell+1)T-1])\leq (1+\mu)^{T-1} V(z[\ell T]).
\end{align}
Now, note that
\begin{align}
V(z[(\ell+1)T])
&=
\begin{bmatrix}
x[(\ell+1)T] \nonumber \\
0
\end{bmatrix}^\top
P\begin{bmatrix}
x[(\ell+1)T] \\
0
\end{bmatrix}\nonumber\\
&=
z[(\ell+1)T-1]^\top F^\top\nonumber \\
&\hspace{.2in}\times\begin{bmatrix}
I & 0 \\
0 & 0
\end{bmatrix}
P\begin{bmatrix}
I & 0 \\
0 & 0
\end{bmatrix} F z[(\ell+1)T-1]\nonumber\\
&\leq \delta z[(\ell+1)T-1]^\top
P z[(\ell+1)T-1]\nonumber\\
&=\delta V(z[(\ell+1)T-1]).\label{eqn:proof:2}
\end{align}
Combining~\eqref{eqn:proof:1} and~\eqref{eqn:proof:2} results in
$
V(z[(\ell+1)T])\leq \delta(1+\mu)^{T-1} V(z[\ell T]).
$
This shows that $\lim_{\ell\rightarrow\infty} V(z[\ell T])=0$. Following~\eqref{eqn:proof:0}, it can be seen that $\lim_{k\rightarrow\infty} V(z[k])=0$. This proves the stability of the system.

\section{Proof of Proposition~\ref{prop:stars}} \label{proof:prop:stars}
\revise{We state the proof for two cases where $n_c=n_x$ and $n_c>n_x$. }\revise{Case I ($n_c=n_x$): } \revise{Since} $(A,B)$ is stabilizable and $(A,C)$ is detectable, a Luenberger observer \revise{exists such that} conditions~\eqref{eqn:tho:1:cond1} and~\eqref{eqn:tho:1:cond2} \revise{are satisfied}.
\revise{Case II ($n_c>n_x$):} Unobservable and uncontrollable states can be \revise{incorporated into} the controller in addition to the Luenberger observer \revise{so that \eqref{eqn:tho:1:cond1} and~\eqref{eqn:tho:1:cond2} are satisfied}.

\section{Proof of Lemma~\ref{lemma:bound}} \label{proof:lemma:bound}
The block $\tilde{\mathbf{F}}(\nu)$ can be factored as
$\tilde{\mathbf{F}}(\nu) = \mathbf{P}(\nu) \begin{pmatrix} 0 & I_n  \end{pmatrix}^\top \mathbf{R}(\nu).$ Then, by properties of the Schur complement, $\mathbf{S}(\nu) \succeq 0$ if and only if
\begin{equation}\label{Synthesis32}
\delta \mathbf{P}(\nu) - \tilde{\mathbf{F}}(\nu)^\top  \mathbf{P}(\nu)^{-1} \tilde{\mathbf{F}}(\nu) = \delta \mathbf{P}(\nu) - \mathbf{R}(\nu)^\top  X \mathbf{R}(\nu) \succeq 0;
\end{equation}
therefore, using the Schur complement again, $\mathbf{S}(\nu) \succeq 0$ if and only if
\begin{equation}\label{Synthesis33}
\begin{bmatrix}  \delta \mathbf{P}(\nu) & \hspace{2mm} \mathbf{R}(\nu)^\top  \\ \mathbf{R}(\nu) & X^{-1} \end{bmatrix} \succeq 0.
\end{equation}
Because $(X^{-1/2} - X^{1/2})^\top (X^{-1/2} - X^{1/2}) \succeq 0$ for any $X$, we have that $X^{-1} \succeq 2I_n - X$. It follows that
\begin{equation*}
\begin{bmatrix}  \delta \mathbf{P}(\nu) & \hspace{2mm} \mathbf{R}(\nu)^\top  \\ \mathbf{R}(\nu) & X^{-1} \end{bmatrix} \succeq \tilde{\mathbf{L}}(\nu) = \begin{bmatrix}  \delta \mathbf{P}(\nu) & \hspace{2mm} \mathbf{R}(\nu)^\top  \\ \mathbf{R}(\nu) & 2I_n - X \end{bmatrix};
\end{equation*}
therefore, $\tilde{\mathbf{L}}(\nu) \succeq 0 \Rightarrow \mathbf{S}(\nu) \succeq 0$.

\section{Proof of Lemma~\ref{lemma:congruence}}\label{proof:lemma:congruence}
Let $\mathbf{P}(\nu)$ be positive definite; then, by properties of the Schur complement, $Y \revise{\,\succ\,}0$ and $X - Y^{-1} \revise{\,\succ\,}0$, and because $YX + VU^\top  = I$ by construction (see  \eqref{Synthesis2}), $VU^\top  = I - YX  \revise{\,\prec\,} 0$, i.e., the matrix $VU^\top $ is nonsingular. Therefore, if $\mathbf{P}(\nu) \revise{\,\succ\,}0$, it is always possible to find nonsingular $U$ and $V$ satisfying $YX + VU^\top  = I$. \revise{The existence of a nonsingular} $V$ implies that $\mathcal{Y}$\revise{, introduced in~\eqref{Synthesis2},} and $\mathcal{T}\revise{:=\text{diag}[\mathcal{Y},\mathcal{Y}]}$ are invertible. Moreover, nonsingular $U$ and $V$ imply that the matrices:
\[
\begin{pmatrix} U & XB \\ 0 & I_{n_u}  \end{pmatrix} \text{ and } \begin{pmatrix} V^\top  & 0 \\ CY & I_{n_y} \end{pmatrix},
\]
are invertible. Therefore, the controller matrices:
\begin{align}\label{CONTROLLER}
\begin{pmatrix} A_c & B_c\\ C_c & D_c \end{pmatrix} = &\begin{pmatrix} U & XB \\ 0 & I_{n_u}  \end{pmatrix}^{-1} \begin{pmatrix} K_1- XAY & K_2 \\ K_3 & K_4 \end{pmatrix}\\ \notag &\times \begin{pmatrix} V^\top  & 0 \\ CY & I_{n_y} \end{pmatrix}^{-1},
\end{align}
are the unique solution of the matrix equation \eqref{change_of_coordinates}.

\section{Proof of Lemma~\ref{synthesis_lemma1}}\label{proof:synthesis_lemma1}
Assume that $\nu$ is such that  $\mathbf{P}(\nu)  \revise{\,\succ\,} 0, \text{ }\mathbf{L}(\nu)  \revise{\,\succeq\,} 0$, and $\tilde{\mathbf{L}}(\nu)  \revise{\,\succeq\,} 0$. Then, by Lemma \ref{lemma:congruence}, $\mathcal{Y}$ and $\mathcal{T}$ are square and nonsingular and thus the transformations $P \rightarrow \mathcal{Y}^\top  P\mathcal{Y} = \mathbf{P}(\nu)$, $\mathcal{L} \rightarrow \mathcal{T}^\top \mathcal{L} \mathcal{T} = \mathbf{L}(\nu)$, and $\tilde{\mathcal{L}} \rightarrow \mathcal{T}^\top \tilde{\mathcal{L}} \mathcal{T} = \mathbf{S}(\nu)$  are congruent. By Lemma \ref{lemma:bound}, $\tilde{\mathbf{L}}(\nu) \revise{\,\succeq\,} 0 \Rightarrow \mathbf{S}(\nu)  \revise{\,\succeq\,}0$. It follows that ($\mathbf{P}(\nu)  \revise{\,\succ\,} 0, \mathbf{L}(\nu)  \revise{\,\succeq\,} 0$, $\tilde{\mathbf{L}}(\nu)  \revise{\,\succeq\,} 0$) \revise{implies} ($P \revise{\,\succ\,}0,\mathcal{L}  \revise{\,\succeq\,} 0,\tilde{\mathcal{L}}  \revise{\,\succeq\,} 0$) because $(\mathbf{P}(\nu),\mathbf{L}(\nu),\mathbf{S}(\nu))$ have the same signature as ($P,\mathcal{L},\tilde{\mathcal{L}}$), respectively, and $\tilde{\mathbf{L}}(\nu) \revise{\,\succeq\,} 0 \Rightarrow \mathbf{S}(\nu)  \revise{\,\succeq\,} 0$. Because $\mathbf{P}(\nu)  \revise{\,\succ\,} 0$, the matrices $U$ and $V$ are nonsingular. This implies that the change of variables in \eqref{change_of_coordinates} and $\mathcal{T}$ are invertible and lead to unique $(P,A_c,B_c,C_c,D_c)$ by inverting \eqref{Synthesis5} and \eqref{change_of_coordinates}.

\section{Proof of Theorem~\ref{synthesis_theorem}}\label{proof:synthesis_theorem}
If $\nu$ satisfies $\mathbf{P}(\nu)  \revise{\,\succ\,} 0$, $\mathbf{L}(\nu)  \revise{\,\succeq\,} 0$, and $\tilde{\mathbf{L}}(\nu) \revise{\,\succeq\,} 0$, by Lemma~\ref{synthesis_lemma1}, the change of variables in \eqref{change_of_coordinates} and matrix $\mathcal{Y}$ in \eqref{Synthesis2} are invertible and the controller matrices \eqref{CONTROLLER} and $P$ obtained by inverting \eqref{Synthesis5} and \eqref{change_of_coordinates} are unique and satisfy inequalities \eqref{eqn:tho:1:cond1}-\eqref{eqn:tho:1:cond3} with $\epsilon = \lambda_{\min}(P)$. Hence, because \eqref{eqn:tho:1:cond4} is satisfied by assumption, by Theorem \ref{tho:1}, the controller matrices in \eqref{CONTROLLER} render the closed-loop dynamics \eqref{eqn:system}-\eqref{eqn:controller} asymptotically stable.

\section{Proof of Theorem~\ref{tho:3}}\label{proof:tho:3}
For the sake \revise{of} brevity $F(\mathcal{P},\bar{\mathcal{C}})$ is denoted by $\bar{F}$. First, note that
\begin{align*}
z[\ell T+i]=\bar{F}^i z[\ell T]+\sum_{j=0}^{i-1}\bar{F}^{i-j-1}G\theta_{\ell T+j},
\end{align*}
for all $i\in\{1,\dots,T-1\}$, where $\theta_k\revise{:=}\bar{y}[k]-y[k]$ and
\begin{align*}
G=
\begin{bmatrix}
A+B\bar{D}_c \\
\bar{B}_c
\end{bmatrix}.
\end{align*}
For $i=T$, it can be seen that
\begin{align*}
z[\ell T+i]=\begin{bmatrix}
I & 0 \\
0 & 0
\end{bmatrix}\bigg(\bar{F}z[\ell T+i-1]+G\theta_{\ell T+i-1}\bigg).
\end{align*}
Combining these update rules shows that
\begin{align} \label{eqn:update_lifted_z}
z[(\ell+1) T]=H\bar{F}^T z[\ell T]+w_\ell,
\end{align}
where
\begin{align*}
H:=\begin{bmatrix}
I & 0 \\
0 & 0
\end{bmatrix},\quad
w_\ell:=\sum_{j=0}^{T-1} \begin{bmatrix}
I & 0 \\
0 & 0
\end{bmatrix}\bar{F}^{T-j-1}G\theta_{\ell T+j}.
\end{align*}
Note that
\begin{align*}
z[(\ell+1) T]^\top Pz[(\ell+1) T]
=&z[\ell T]^\top \bar{F}^{\top T} H^\top P H \bar{F}^T
z[\ell T]\\
&+2w_\ell^\top P H \bar{F}^Tz[\ell T]+w_\ell^\top P w_\ell\\
\leq &\delta z[\ell T]^\top \bar{F}^{\top T-1} P  \bar{F}^{T-1}
z[\ell T]\\
&+2w_\ell^\top P H \bar{F}^Tz[\ell T]+w_\ell^\top P w_\ell\\
\leq &\delta(1+\mu)^{T-1} z[\ell T]^\top P z[\ell T]\\
&+2w_\ell^\top P H \bar{F}^Tz[\ell T]+w_\ell^\top P w_\ell,
\end{align*}
where the inequalities follow from $\bar{F}^{\top } H^\top P H \bar{F}\leq \delta P$ and $\bar{F}^{\top }P \bar{F}\leq (1+\mu) P$. Thus,
\begin{align*}
z[(\ell+1) T]^\top &Pz[(\ell+1) T]
-z[\ell T]^\top P z[\ell T]\\
\leq &(\delta(1+\mu)^{T-1}-1) z[\ell T]^\top P z[\ell T]\\
&+2w_\ell^\top P H \bar{F}^Tz[\ell T]+w_\ell^\top P w_\ell\\
\leq &(\delta(1+\mu)^{T-1}-1) \|P^{1/2}z[\ell T]\|_2^2\\
&+2\|P^{-1/2}\bar{F}^{\top T} H^\top P\|_2 W\|P^{1/2}z[\ell T]\|_2\\
&+\lambda_{\max}(P) W^2,
\end{align*}
where
\begin{align*}
W
:=&\sup_{\ell} \|w_\ell\|_2\\
=&\sup_{\ell} \bigg\|\sum_{j=0}^{T-1} \begin{bmatrix}
I & 0 \\
0 & 0
\end{bmatrix}\bar{F}^{T-j-1}G\theta_{\ell T+j}\bigg\|_2\\
\leq & \sum_{j=0}^{T-1} \bigg\|\begin{bmatrix}
I & 0 \\
0 & 0
\end{bmatrix}\bar{F}^{T-j-1}G\theta_{\ell T+j}\bigg\|_2\\
\leq &\sum_{j=0}^{T-1} \|\bar{F}\|^{T-j-1}\|G\|\sup_{\ell}\|\theta_{\ell T+j}\|_2\\
\leq &\bigg(\sum_{j=0}^{T-1} \|\bar{F}\|^{T-j-1}\|G\|\bigg)\revise{n_y}2^{-m},
\end{align*}
where the last equality follows from that the quantization error is bounded by $\|\theta_{\ell T+j}\|_2\leq \revise{n_y} 2^{-m}$ if $n$ is selected large enough\revise{, i.e., $n$ is selected such that $\|y[k]\|_\infty< 2^{n-1}$ (see the last paragraph of this proof for ensuring this)}. Let $s_1$ denote the largest real root of the polynomial equation $s\mapsto [(\delta(1+\mu)^{T-1}-1)\lambda_{\max}(P)] s^2+[2\|P^{-1/2}\bar{F}^{\top T} H^\top P\|_2 W]s+ [\lambda_{\max}(P) W^2]$. Define $\mathbb{M}_1:=\{z\,|\,z^\top P z\leq s_1^2\}$. Notice that if $z[\ell T]\notin \mathbb{M}_1$, $z[(\ell+1) T]^\top Pz[(\ell+1) T]-z[\ell T]^\top P z[\ell T]<0$ (and thus $z[(\ell+1) T]\in \mathbb{M}_1$) since $(\delta(1+\mu)^{T-1}-1\revise{)}<0$.  However, if $z[\ell T]\in \mathbb{M}_1$, it can be deduce that
$z[(\ell+1) T]^\top Pz[(\ell+1) T]
\leq s_2
:=\delta(1+\mu)^{T-1} s_1^2+\lambda_{\max}(P) W^2
+2\|P^{-1/2}\bar{F}^{\top T} H^\top P\|_2 Ws_1.$
Define $\mathbb{M}_2:=\{z\,|\,z^\top P z\leq s_2^2\}$. Evidently, if $z[\ell T]\in \mathbb{M}_1$, then $z[(\ell+1) T]\in \mathbb{M}_2$. Combining these results, it can be seen that $\mathbb{M}_1\cup \mathbb{M}_2$ is an invariant set for~\eqref{eqn:update_lifted_z}. This is because two distinct cases can happen if $z[\ell T]\revise{\,\in\,}\mathbb{M}_1\cup \mathbb{M}_2$; either $z[\ell T]\in \mathbb{M}_2\setminus \mathbb{M}_1$ or $z[\ell T]\in \mathbb{M}_1$ must happen. If $z[\ell T]\in \mathbb{M}_2\setminus \mathbb{M}_1$, it means that $z[\ell T]\notin \mathbb{M}_1$, thus $z[(\ell+1) T]\in \mathbb{M}_1\subseteq \mathbb{M}_1\cup \mathbb{M}_2$. On the other hand, if $z[\ell T]\in \mathbb{M}_1$, then $z[(\ell+1) T]\in \mathbb{M}_2\subseteq \mathbb{M}_1\cup \mathbb{M}_2$. Furthermore, $\mathbb{M}_1\cup \mathbb{M}_2$ is an invariant set for the dynamical system in feedback loop with~\eqref{eqn:controller_quantized}. This is because, for all $k\in\{\ell T,\dots,(\ell+1)T-2\}$, $z[k+1]^\top P z[k+1]-z[k]^\top P z[k]=V(z[k+1])-V(z[k])\leq \mu z[k]^\top P z[k]<0$. It can be seen that the set $\mathbb{M}_1\cup \mathbb{M}_2$ is attractive for~\eqref{eqn:update_lifted_z}. This because if $z[\ell T]\notin \mathbb{M}_1\cup \mathbb{M}_2$, it must also $z[\ell T]\notin \mathbb{M}_1$. Therefore, $z[(\ell+1) T]^\top Pz[(\ell+1) T]-z[\ell T]^\top P z[\ell T]<0$. All these results in the fact that $\lim_{k\rightarrow\infty} \dist (z[k],\mathbb{M}_1\cup \mathbb{M}_2)=0$. Now, note that $\max(s_1,s_2)=\mathcal{O}(2^{-m})$. Therefore, there exists $\varrho>0$ such that $\lim_{k\rightarrow\infty} \dist (z[k],\mathcal{B}(\varrho 2^{-m}))=0$.

It only remains to find the bound on $n$. Note that the largest value of the Lyapunov function is
\begin{align*}
z[0]^\top Pz[0]=c:=x_0^\top \begin{bmatrix}
I & 0 \\
0 & 0
\end{bmatrix}P\begin{bmatrix}
I & 0 \\
0 & 0
\end{bmatrix}x_0.
\end{align*}
This implies that $z[k]^\top Pz[k]\leq c$. Therefore, $x[k]^\top x[k]+x_c[k]^\top x_c[k]=z[k]^\top z[k]\leq c/\lambda_{\min}(P):=c/\epsilon$. As a result, $x[k]^\top x[k]\leq c/\epsilon$. Noting that $y[k]=Cx[k]$, it can be deduced that $y[k]^\top y[k]\leq \lambda_{\max}(C^\top C)c/\epsilon$. Finally, because of the relationship between norms, it can be seen that $\|y[k]\|_{\infty}\leq \lambda_{\max}(C^\top C)c/\epsilon$. This implies that $2^{n-1}> \lambda_{\max}(C^\top C)c/\epsilon$ and thus $n>\log_2(\lambda_{\max}(C^\top C)c/\epsilon)+1$.

\section{Proof of Lemma~\ref{lemma:2}} \label{proof:lemma:2}
Noting that the controller resets every $T$ steps, we only need to prove this result for $k\in\{0,\dots,T-1\}$. At $k=0$, $x_c[k]=0$ and thus $x_c[1]\in\mathbb{Q}(n_y+2n,2m)$ because the entries of $\bar{B}_c\bar{y}[k]$ belong to $\mathbb{Q}(n_y+2n,2m)$. For all $k\in\{1,\dots,T-1\}$, if the entries of $x_c[k]$ belong to $\mathbb{Q}(n',m')$, the entries of $\bar{A}_cx_c[k]$ (at worst case) belong to $\mathbb{Q}(n_c+n+n',m+m')$ and the entries of $\bar{B}_c\bar{y}[k]$ belong to $\mathbb{Q}(n_y+2n,2m)$, therefore the entries of $x[k+1]=\bar{A}_cx_c[k]+\bar{B}_c\bar{y}[k]$ must belong to $\mathbb{Q}(n_c+n+n'+1,m+m')$ because $n_y\leq n_c$, $n\leq n'$, and $m\leq m'$. Furthermore, $u[k]$ must belong to $\mathbb{Q}(n_c+n+n'+1,m+m')$. This proves that the entries of $x_c[k]$ and $u[k]$ must, respectively, belong to $\mathbb{Q}((n_c+n+1)(k\modd T-1)+n_y+2n,m(k\modd T-1)+2m)$ and $\mathbb{Q}((n_c+n+1)(k\modd T)+n_y+2n,m(k\modd T)+2m)$.

\bibliographystyle{ieeetr}
\bibliography{ref}

\end{document}